%% file: main.tex
\documentclass{article}

\usepackage[main, final]{neurips_2026}

\input{common}

\usepackage[utf8]{inputenc} 
\usepackage[T1]{fontenc}    
\usepackage{url}            
\usepackage{booktabs}       
\usepackage{amsfonts}       
\usepackage{nicefrac}       
\usepackage{microtype}      
\usepackage{xcolor}         

\title{Learning inexact alternating minimization \\ for convex problems with quadratic coupling}

\author{%
  Paul Häusner$^{1}$
  \qquad
  Jevgenija Rudzusika$^{2}$
  \qquad
  Jens Sjölund$^1$
  \qquad
  Ozan Öktem$^2$ 
  \\ \\
  $^1$Department of Information Technology, Uppsala University, Sweden \\
  $^2$Department of Mathematics, KTH Royal Institute of Technology, Stockholm, Sweden \\
  \texttt{paul.hausner@it.uu.se, jevaks@kth.se, jens.sjolund@it.uu.se, ozan@kth.se}
}

\begin{document}

\maketitle

\begin{abstract}
Alternating minimization is a standard tool for optimization problems with two coupled variable blocks.
Its cost is dominated by the inner subproblem solvers, which rarely have a closed-form solution.
In this paper, we propose to learn an inexact alternating scheme in which neural networks approximate the subproblem solution at each iteration.
For convex problems with quadratic couplings between the variable blocks that satisfy a quadratic growth condition, we derive a worst-case bound on the optimality gap that exposes the per-iteration error that each network is trained to minimize in a greedy fashion.
Combining the bound with a PAC-Bayes argument, we obtain a fixed-horizon, high-probability bound on the expected optimality gap over the distribution of problem instances.
We apply the learned scheme to robust matrix completion and low-dose CT reconstruction with dictionary regularization where we achieve a $5\times$--$13\times$ speedup over the inexact first-order baseline at matched accuracy.
A brief end-to-end fine-tuning stage extends this up to $30\times$ and outperforms other learned schemes.
The code is available at \url{https://github.com/paulhausner/learned-alternating-minimization}.
\end{abstract}

\section{Introduction}
We consider large-scale convex optimization problems where one seeks to minimize a convex objective with two variable blocks 
$(x,z) \mapsto F(x,z)$ of the form\looseness=-1
\begin{equation}
\label{eq:Objective}
 F(x, z) = f(x) + \frac{1}{2}\|A x + B z \|^2 + g(z)
  \quad\text{for $(x,z) \in \mathcal{X} \times \mathcal{Z}$.}
\end{equation}
Here, $\mathcal{X}$ and $\mathcal{Z}$ are Hilbert spaces.
Furthermore, both $f \colon \mathcal{X} \to \Real$ and $g \colon \mathcal{Z} \to \Real$ are 
proper, convex, lower semi-continuous functions. 
The two variable blocks $x$ and $z$ are coupled through a quadratic function that is defined in terms of 
bounded linear operators $A$ and $B$ that map $\mathcal{X}$ and $\mathcal{Z}$ to a common Hilbert space with inner product $\langle \cdot, \cdot\rangle$ and corresponding induced norm $\| \cdot \|$.
Such problems appear in, e.g., inverse problems \citep{cheng2020image, chambolle2015remark}, control \citep{khan2025half}, and naturally arise in half-quadratic splitting formulations \citep{geman1995nonlinear}.

Directly optimizing the objective function~\cref{eq:Objective} is challenging due to the quadratic coupling between the two variable blocks that prevents efficient joint optimization.
The following alternating minimization scheme \citep{tseng2001convergence} instead decouples the problem into two subproblems each involving only a single variable block that can be solved more efficiently using first-order methods: 
\begin{align}
    \label{scheme1}
    x^{k+1} &:= \argmin_{x \in \mathcal{X}} F(x, z^k) &
    z^{k+1} := \argmin_{z \in \mathcal{Z}} F(x^{k+1}, z).
\end{align}
This alternating scheme is well studied and achieves a sublinear convergence with rate $\mathcal{O}(1/K)$ \citep{beck2017first}
or under additional assumptions on the strong convexity or regularity conditions \citep{liao2024error},
a linear convergence $\mathcal{O}(\rho^K)$ with rate $\rho < 1$ \citep{beck2013convergence}.\looseness=-1

The convergence speed described in the previous paragraph refers to the number of \emph{outer iterations}~$k$.
However, the two subproblems in an outer iteration \cref{scheme1} 
do not necessarily have a closed-form solution and are instead often solved with an iterative scheme themselves.
We refer to these as \emph{inner iterations}.
Using an iterative solver for the inner subproblems has both practical and theoretical disadvantages: it is computationally costly (resulting in long wall-clock solving times) and limits the usability of many theoretical results, which often assume the subproblems are solved exactly.
Our main goal is to replace the updates in \cref{scheme1} with an inexact solver that can be implemented, for example, using a neural network.
This allows us to compute the next iterate significantly faster and avoid the computational bottleneck.
In particular, our contributions can be summarized as follows:\looseness=-1
\begin{itemize}
    \item In \cref{sec:theory}, we derive a worst-case error bound of the inexact alternating scheme (\cref{theorem:linearbound}) under the quadratic growth assumption and prove linear and sublinear convergence under suitable conditions on the decay of the error terms (\cref{corollary2}). 
    \item Based on the inexact alternating scheme, we take a probabilistic viewpoint and derive a novel learning-based approach to accelerate inexact alternating minimization by predicting the solution to the inner subproblems. We propose a greedy training strategy that optimizes the networks to solve each subproblem in isolation (\cref{sec:greedytr}).
    \item We derive a fixed-horizon PAC-Bayes bound on the expected optimality gap of the learned scheme under the posterior over its weights (\cref{thm:pacbayes}). With a random-walk prior over the per-iteration weights, the complexity term reduces to the parameter drift across iterations.
\end{itemize}
In numerical experiments (\cref{sec:ctresults,sec:matrixresults}), we apply both inexact and learned methods to solve CT reconstruction with dictionary-based regularization \citep{rudzusika2022deep} and robust matrix completion \citep{tao2011recovering}.
By predicting the solution to the subproblems with the greedily-trained networks we accelerate the inexact alternating scheme by an order of magnitude at matched accuracy.

\section{Related work}
\label{sec:relatedwork}

\subsection{Data-driven optimization}
Data-driven optimization methods are justified in settings where structurally similar problems need to be solved repeatedly \citep{chen2022learning}.
This allows the upfront cost of training to be amortized by solving many subsequent problems of the same type.
Learned methods show strong empirical performance \citep{li2016learning} but obtaining guarantees can be challenging \citep{amos2023tutorial}.\looseness=-1

\myparagraph{Learning initializations}
One natural approach is to learn the initial guess for an optimization algorithm. 
\citet{sambharya2024learning} train neural networks end-to-end to warm start fixed point optimization algorithms, and \citet{sjolund2022graph} predict the initialization of an alternating primal-dual scheme. 
Along similar lines, \citet{schmidtobreick2025warm} learn to warm start dual active-set methods.
These approaches inherit the convergence properties of existing methods.
For alternating schemes a good starting guess is readily available as the solution to the subproblem from the previous iteration.\looseness=-1

\myparagraph{Learning iterative updates}
To derive convergence guarantees for learned updates, data-driven optimization methods restrict the allowed predictions.
The safeguarding scheme proposed by \citet{heaton2023safeguarded} obtains convergence of the learned iterates using the theory of fixed-point operators, which results in asymptotic convergence.
Similarly, \citet{premont2022simple} prove the asymptotic convergence of their proposed safeguarding mechanism based on 
choosing safe updates of the iterates.
\citet{liu2023towards} formalize and provide necessary conditions on the learned iterates of first-order methods 
to ensure asymptotic convergence that can be implemented for practical problems, for example, by learning preconditioners 
\citep{hausner2024neural, doerks2025learning}.
Along similar lines to our training strategy, \citet{fahy2024greedy} propose to learn a preconditioner in a greedy fashion.

\citet{banert2020data} learn the parameters of optimization algorithms, enforcing certain conditions on the learned parameters to retain the asymptotic convergence behavior.
\citet{sambharya2024learning2} learn hyper-parameters of optimization algorithms but require reversion to a safe parameter region to inherit asymptotic convergence.
Further, \citet{banert2024accelerated} learn a deviation that restricts the updates to some a priori known set of provably convergent updates.
Similar to our approach, \citet{tan2023data} provide a worst-case analysis of their learned mirror descent framework, but the obtained bound does not guarantee asymptotic convergence in general.
For QP solvers, \citet{gao2025learning} propose a sequence-based approach to learn the subproblem solution.
We instead formulate a learning problem for each iteration in isolation and analyze the combined convergence of separate networks.
By considering the continuous-time dynamics of optimization algorithms, \citet{xie2024ode} derive a linear convergence rate of learned updates.
Further, \citet{martin2025learning} characterize learned optimizers with a linear convergence rate.
More recently, \citet{sambharya2024data} and \citet{sucker2024learning} propose to use data-driven probabilistic guarantees both to estimate the performance of classical optimizers on parametric problems and to directly train learned optimizers on the PAC-Bayes guarantees.\looseness=-1

\myparagraph{Learned image reconstruction}
Instead of first deriving an optimization problem that is then solved, several approaches directly integrate the learning and optimization by incorporating the regularizer implicitly into the learned optimizer, which limits both interpretability and theoretical analysis \citep{10004773}.
The seminal work by \citet{gregor2010learning} learns sparse coding and \citet{adler2018learned} perform reconstruction with a neural network obtained by unrolling the primal-dual hybrid gradient method for minimizing the $\ell_2$ distance.
For alternating schemes, \citet{ding2023learned} learn the regularization term instead of a solution to an optimization problem similar to \citet{xin2022learned} where the method is applied to half-quadratic splitting.
Further, \citet{fahes2021unrolling} unroll a linearized alternating scheme to learn the algorithm's hyperparameters and variables. \looseness=-1

\subsection{Inexact optimization}
Inexact optimization methods have a long history of explicitly modeling suboptimal solutions to the subproblems while maintaining asymptotic convergence or proving worst-case bounds under the common assumption of error summability \citep{guler1992new}.
In the context of alternating minimization, \citet{robini2016inexact} solve related half-quadratic splitting problems inexactly under KL regularity conditions.
The related primal-dual ADMM \citep{eckstein1992douglas} has also been studied extensively in the inexact setting \citep{ouyang2015accelerated, xie2017inexact}.
However, we restrict ourselves to a class of alternating optimization schemes that naturally arise in inverse problems such as sparse signal recovery \citep{Wang:2024} and automatic control \citep{Aravkin2014}.\looseness=-1

\section{Inexact scheme}
\label{sec:theory}

In this section, we propose an inexact scheme to solve alternating minimization problems with quadratic coupling terms between 
the variable blocks as in \cref{eq:Objective}.
Specifically, in each iteration $k$, we allow for inexact solutions $\tilde{x}^{k+1}$ and $\tilde{z}^{k+1}$ to the subproblems in the alternating scheme \cref{scheme1}:
\begin{align}
    \label{scheme2}
    \tilde{x}^{k+1} &\approx \argmin_{x \in \mathcal{X}} F(x, \tilde{z}^k),  &  \tilde{z}^{k+1} \approx \argmin_{z \in \mathcal{Z}} F(\tilde{x}^{k+1}, z).
\end{align}
To measure the performance of the inexact scheme and derive convergence results, we denote the function gap introduced by the inexact scheme in each step respectively as:
\begin{align}
    \label{eq:inexact-error}
    \xi^{k+1} &:= F(\tilde{x}^{k+1}, \tilde{z}^{k}) - F(x^{k+1}, \tilde{z}^{k}) & 
    \zeta^{k+1} := F(\tilde{x}^{k+1}, \tilde{z}^{k+1}) - F(\tilde{x}^{k+1}, z^{k+1}),
\end{align}
where $x^{k+1} := \argmin_{x \in \mathcal{X}} F(x, \tilde{z}^k)$ and $z^{k+1} := \argmin_{z \in \mathcal{Z}} F(\tilde{x}^{k+1}, z)$ are the exact minimizers of the respective subproblems.
Both quantities are non-negative by construction.
To simplify the analysis of the inexact scheme, we also introduce the combined error $\Delta^{k+1} := \xi^{k+1} + \zeta^{k+1}$.

Next, we derive a worst-case bound of the optimality gap based on the iterates from the alternating scheme and the sequence of combined errors $\{ \Delta^k \}_{k=1}^K$ introduced by the inexact solvers.
The key ingredient in the proof of the worst-case bounds from \cref{theorem:linearbound} and convergence is the following lemma, which shows a one-step inequality for the iterates of the inexact alternating scheme.
\begin{lemma}
    \label{lemma0p}
    Assume that $F$ is a convex objective of the form \cref{eq:Objective} and that the iterates $\{ \tilde{x}^k, \tilde{z}^k \}_{k=0}^K \subseteq \mathcal{X} \times \mathcal{Z}$ are obtained using the inexact scheme \cref{scheme2}.
    Then for any iteration~$k \geq 0$ it holds for any $(x, z) \in \mathcal{X} \times \mathcal{Z}$ that the new function value $\tilde{F}^{k+1} := F(\tilde{x}^{k+1}, \tilde{z}^{k+1} )$ satisfies\looseness=-1
    \begin{equation}
        \label{eq:lyapanov}
         \tilde{F}^{k+1} + \frac{1}{2(1+\alpha)}\|B(z - \tilde{z}^{k+1})\|^2 \leq F(x, z )
        + \frac{1+\alpha}{2}\|B(z - \tilde{z}^{k})\|^2 + \frac{1+\alpha}{\alpha}\, \Delta^{k+1}
    \end{equation}
    where $\alpha > 0$ is an arbitrarily chosen constant.
\end{lemma}

To prove convergence of the inexact scheme we utilize the following assumption which states that the function is bounded from below by a quadratic form in the $B$-norm distance to the optimal set.
Let $Z^* := \argmin_{(x, z) \in \mathcal{X} \times \mathcal{Z}} F(x, z)$ denote the set of minimizers and $F^*$ the optimal value.\looseness=-1

\begin{assumption}[Quadratic growth]
    \label{assumption:quadraticgrowth}
    We assume that the function gap between the inexact iterates $\{\tilde{x}^k, \tilde{z}^k \}$ and the optimal set $Z^*$ can be bounded in $B$-norm distance with parameter $\mu > 0$ as:\looseness=-1
    \begin{equation}
        \frac{\mu}{2} \inf_{(x', z') \in Z^*} \| B(z' - \tilde{z}^k)\|^2 \leq F(\tilde{x}^k, \tilde{z}^k) - F^*.
    \end{equation}
\end{assumption}
This assumption is implied if the function $F(x, z)$ is jointly strongly convex in both arguments but belongs to a weaker class of assumptions since the inequality only has to hold around the optimal set \citep{liao2024error}.
Now, we present the following theorem, which establishes a worst-case bound on the objective in terms of the initial iterates and the errors introduced by the inexact scheme.
\begin{theorem}[Worst-case bound for the inexact scheme]
    \label{theorem:linearbound}
    Let the objective function $F$ in \cref{eq:Objective} be convex with $F^*$ as its global minimum value.
    The iterates $\{ \tilde{x}^k, \tilde{z}^k \}_{k=0}^K \subseteq \mathcal{X} \times \mathcal{Z}$ are obtained from the inexact alternating scheme in \cref{scheme2} with an arbitrary initial $(\tilde{x}^0, \tilde{z}^0) \in \mathcal{X} \times \mathcal{Z}$.
    Under \cref{assumption:quadraticgrowth} with some fixed $\mu$ and choosing the parameter $\alpha$ as in \cref{lemma0p} from the interval
    \begin{equation}
        \label{eq:alphabounds}
        \alpha \in \left( 0, \tfrac{\mu - 2 + \sqrt{\mu^2 + 4}}{2} \right),
    \end{equation}
    and with the initial Lyapunov energy function defined as
    \begin{equation}
        \Phi_0 := \tilde F^0 - F^* + \tfrac{1}{2(1+\alpha)}\inf_{(x', z') \in Z^*} \|B(z' - \tilde z^0)\|^2,
    \end{equation}
    we can bound the function gap for any iteration $K > 0$ as:
    \begin{equation}
        \label{eq:theorem1}
        \tilde{F}^K - F^* \leq \rho^{K} \Phi_0 + \frac{1+\alpha}{\alpha} \sum_{k=1}^K \rho^{K-k} \Delta^k,
    \end{equation}
    where $\rho = {(1+\alpha)^2} / {(1+\mu(1+\alpha))} < 1$ is the contraction rate of the bound.
\end{theorem}

\begin{remark}
    The parameter $\alpha$ in \cref{theorem:linearbound} is a hyper-parameter of the bound that can be used to trade-off between the solver inaccuracy and the convergence rate of the error bound.
    In particular, if the errors $\Delta^k$ are zero, we can choose $\alpha \to 0$ and recover the exact linear convergence rate $\rho \to 1/(1+\mu)$.
    Conversely, increasing $\alpha$ reduces the amplification of the error term $(1+\alpha)/\alpha$ but worsens the contraction rate $\rho$ towards $1$ slowing down the outer convergence.
\end{remark}
The first term in the bound decays geometrically with rate $\rho$, so to establish convergence of the scheme, we are required to analyze how the error $\Delta^{k+1}$ decays.
Taking into account the structure of the accumulated errors, we can make the following key observation before deriving the convergence.
\begin{remark}
    \label{remark:discount}
    In the bound from \cref{theorem:linearbound}, early error terms are discounted exponentially.
    Thus solvers that become more accurate over time can compensate for large errors made early on.
\end{remark}

\begin{corollary}
    \label{corollary2}
    Under the same assumptions as in \cref{theorem:linearbound}, the convergence of the function value is determined by the decay of the combined error $\Delta^k$.
    We distinguish three separate regimes:\looseness=-1
    
        \textbf{i.} \emph{Geometric decay.} If $\Delta^k \leq Cq^k$ for some constants $C > 0$ and $q \in (0, 1)$, we retain linear convergence of the inexact scheme with rate $\mathcal{O}(\max(\rho, q)^K)$.
        
        \textbf{ii.} \emph{Polynomial decay.} If $\Delta^k \to 0$, we achieve convergence to the optimum where the sublinear rate is determined by the speed of the decay. In particular, for polynomial decay $\Delta^k \leq C/k^q$ with $q > 0$, we obtain sublinear convergence with a rate of $\mathcal{O}(1/K^q)$.

        \textbf{iii.} \emph{Eventually bounded errors.} If there exist $\bar k \in \mathbb{N}$ and $\Delta^\text{tail} > 0$ such that $\Delta^k \le \Delta^\text{tail}$ for all $k > \bar k$, the function gap converges to a neighborhood of the optimum: $\limsup_{K \to \infty} \tilde{F}^K - F^* \le \Delta^\text{tail}(1+\alpha) / [\alpha(1-\rho)]$.
\end{corollary}
The final case in \cref{corollary2} is particularly relevant for the learned scheme in the next section.
While we cannot guarantee that the error terms vanish for the neural networks, we provide conditions under which the objective gap for a fixed horizon $K$ is bounded with high probability over the training data.\looseness=-1

\section{Learned scheme}
\label{sec:learnedopt}
In this section, we approximate the solution to each subproblem in \cref{scheme2} with a learned mapping.
We implement each inexact solver using a neural network that takes the current iterate together with the gradient of the corresponding partial objective and outputs the next iterate of the inexact scheme:
\begin{align}
    \label{scheme:learned}
    \tilde{x}^{k+1} = & \; \phi_{\theta^{k+1}} \!\left( \tilde{x}^k,\, \nabla_x F(\tilde{x}^k, \tilde{z}^k) \right)
    & \tilde{z}^{k+1} = \; \psi_{\vartheta^{k+1}} \!\left( \tilde{z}^k,\, \nabla_z F(\tilde{x}^{k+1}, \tilde{z}^k)\right).
\end{align}
Each network thus realizes a learned first-order map for the corresponding subproblem.
For non-smooth problems, the gradient input refers to the gradient of the smooth part while the proximal behavior of the non-smooth part is absorbed into the network architecture.
To ensure the learned iterates are bounded, the last layer of the network projects the data onto a box of predefined dimensions derived from the problem domain by clipping the network outputs $\Pi_{[l, u]}$.
Here, $\theta^{k+1}$ and $\vartheta^{k+1}$ are the neural network parameters of each network at iteration $k+1$.
We collect both into the joint parameter set $\Theta^{k+1} := (\theta^{k+1}, \vartheta^{k+1})$ and write $\Theta^{1:k} := (\Theta^1, \ldots, \Theta^k)$ for the parameters up to iteration $k$.

Convergence results in \cref{sec:theory} apply to a fixed objective.
To learn the solution map, we consider a family of objectives $\{ F_y \}_{y \in \mathcal{Y}}$ indexed by a random variable $y \sim p_y$ with bounded support \citep{amos2023tutorial}.
Each $y$ yields an objective as in \cref{eq:Objective} so iterates $\tilde{x}^k_y, \tilde{z}^k_y$ and errors $\Delta^k_y$ are random variables.
For fixed weights $\Theta^{1:K}$, these are functions of $y$ alone and  expectations $\mathbb{E}[\,\cdot\,]$ are with respect to $p_y$.

For each $y \in \mathcal{Y}$, the bound in \cref{theorem:linearbound} applies to the corresponding instance $F_y$.
To characterize the alternating scheme's behavior over a distribution of instances on which a model can be trained, we take the expectation with respect to $p_y$.
The projection $\Pi_{[l, u]}$ confines the learned iterates to a bounded set for every parameter value.
We assume \cref{assumption:quadraticgrowth} holds on this set with a common constant $\mu$ for all $y$ in the support of $p_y$, so that $\alpha$ and $\rho$ are independent of the instance and network weights.
Together with the bounded support of $p_y$, the following result holds as an analogous version in expectation:\looseness=-1
\begin{equation}
    \label{eq:expectbound}
    \mathbb{E}[\tilde{F}^K_y - F^*_y] \;\leq\; \rho^K\, \mathbb{E} \left[ \Phi_0 \right] + \frac{1+\alpha}{\alpha}\sum_{k=1}^K \rho^{K-k}\, \mathbb{E}[\Delta^k_y].
\end{equation}
If the conditions of \cref{corollary2} hold for the expected errors $\mathbb{E}[\Delta^k_y]$, the same argument and linearity of expectation extend its three convergence regimes to convergence in expectation.

\subsection{Greedy training}
\label{sec:greedytr}
We now use the bound \cref{eq:expectbound} to motivate a greedy training objective.
On the right-hand side, the sum of weighted per-iteration errors is the only term that depends on the network.
Minimizing jointly over all $\Theta^{1:K}$ is computationally expensive.
Each error $\Delta^k_y$ depends on the whole trajectory through the previous iterates so the gradient computation requires backpropagation through $k$ unrolled iterations.
Memory therefore grows linearly in $K$ and gradients through long unrolls can be unstable~\citep{metz2019understanding}.
Instead, we only minimize the errors $\Delta^k_y$ greedily one at a time in their order in the scheme.\looseness=-1

For each outer iteration $k$, with $\Theta^{1:k-1}$ held fixed, we first train the $x$-network $\phi_{\theta^k}$ to minimize $\mathbb{E}[\xi^k_y]$.
With $\theta^k$ also held fixed, we then train the $z$-network $\psi_{\vartheta^k}$ to minimize $\mathbb{E}[\zeta^k_y]$.
Since $\zeta^k_y$ also depends on $\theta^k$ through $\tilde x^k_y$, the two phases do not in general minimize $\mathbb{E}[\Delta^k_y]$ jointly.
We refer to this strategy as \emph{greedy training}: each network is trained on its own subproblem, with all previously trained networks held fixed and no anticipation of downstream iterations, so it does not minimize the bound in general.\looseness=-1

Minimizing the population risk requires access to the unknown distribution $p_y$.
We replace the expectation by an empirical mean over i.i.d.\ samples $\mathcal{D} \!=\!\{y_i\}_{i=1}^N$ drawn from the distribution $p_y$.
Since $\Delta^k_y = \xi^k_y + \zeta^k_y$, the training phases minimize, one after the other, the two terms of the empirical risk\looseness=-1
\begin{equation}
    \label{eq:empiricalrisk}
    \hat\Delta^k(\Theta^{1:k}; \mathcal{D}) := \tfrac{1}{N} \sum_{i=1}^N \Delta^k_{y_i}(\Theta^{1:k}) = 
    \tfrac{1}{N}\!\sum_{i=1}^N \xi^k_{y_i}(\Theta^{1:k}) + \tfrac{1}{N}\!\sum_{i=1}^N \zeta^k_{y_i}(\Theta^{1:k}).
\end{equation}
The error terms we minimize are defined in \cref{eq:inexact-error} and also depend on their respective reference values, which are, however, independent of the parameters updated in that phase.
The training loss is therefore equivalent to the objective evaluated at the predicted iterate itself.
This yields an unsupervised training scheme that does not require ground-truth subproblem minimizers as targets.\looseness=-1

\input{algorithms/training}
\Cref{alg:layerwisetraining} formalizes the greedy training procedure described above.
We warm-start each new pair of networks from the weights trained at the previous iteration.
As the iterates approach the minimizer, consecutive subproblems become increasingly similar, so the weights from the previous iteration provide a strong initialization for iteration $k$.
The warm-start lets us cut the training budget at later iterations, where each pair of networks only requires fine-tuning rather than training from scratch.

A natural alternative to the proposed greedy approach is end-to-end training, in which all $K$ networks are optimized jointly through the full alternating scheme.
With end-to-end training, however, the networks are no longer trained to solve the per-iteration subproblems, so the training objective no longer matches the per-iteration errors in the bound and the bound no longer motivates the training.
In addition, greedy training avoids backpropagation through the full $K$-iteration scheme, which reduces memory consumption and accelerates training \citep{metz2019understanding}.

\subsection{Generalization bound}
\label{sec:generalization}
We now bound the expected optimality gap of the learned scheme on unseen instances at a fixed horizon $K$.
Greedy training fits the networks to the empirical risk $\hat\Delta^k$, while the bound in expectation in \cref{eq:expectbound} depends on the population risk $\mathbb{E}[\Delta^k_y]$.
A PAC-Bayes bound \citep{alquier2024user, catoni2007, mcallester1999pac} closes this gap for a randomized predictor drawn from a posterior $\mathbb Q$, at the cost of the divergence $\mathrm{KL}(\mathbb Q \| \mathbb P)$ to a prior $\mathbb P$ fixed before seeing the data \citep{parrado12a}.\looseness=-1


Training on the $N$ training samples $\mathcal{D}$ produces a trajectory of weights $\hat\Theta^1, \ldots, \hat\Theta^K$, one network pair per iteration.
The bound concerns a randomized version of this trajectory: Gaussian noise of variance $\varsigma^2$ is added to each increment $\hat\Theta^k - \hat\Theta^{k-1}$, so that a random walk runs alongside the trained weights.
The posterior $\mathbb Q$ is the distribution of this walk given $\mathcal{D}$.
The prior $\mathbb P$ is the same walk with zero-mean increments from the data-independent initialization $\Theta^0 =: \hat\Theta^0$, which extends the initialization-centered Gaussian prior \citep{dziugaite2017computing} to a sequence of networks.
Since the two walks have the same noise and differ only in their mean increments, their divergence is the cumulative parameter drift $\mathrm{KL} := \mathrm{KL}(\mathbb Q \| \mathbb P) = \tfrac{1}{2\varsigma^2} \sum_{k=1}^K \|\hat\Theta^k - \hat\Theta^{k-1}\|^2$.\looseness=-1

\begin{theorem}[Fixed-horizon generalization bound]
    \label{thm:pacbayes}
    Let $\mathbb P$ and $\mathbb Q$ be the prior and posterior defined above with $\varsigma > 0$ fixed before training.
    Pick a fixed $\alpha$ that satisfies \cref{eq:alphabounds} with the common constant $\mu$, and let $\rho$ be as in \cref{theorem:linearbound}.
    Suppose that $0 \leq \Delta^k_y(\Theta^{1:k}) \leq M$ for all $y$ in the support of $p_y$, all $\Theta^{1:k}$, and $k = 1, \ldots, K$.
    Fix $\beta \in (0, 1)$ and an integer $J \geq 0$, and let $C := \log((2J+1)/\beta)$ and $\lambda^\star := \tfrac{2(1-\rho)}{M}\sqrt{2N(\mathrm{KL}+C)}$.
    Then, with probability at least $1 - \beta$ over the $N$ i.i.d.\ training samples, the following holds: if $2^{-J} \leq \lambda^\star \leq 2^J$, then
    \begin{equation}
      \label{eq:probconv}
      \mathbb E_\mathbb Q \mathbb E_y[\tilde F^K_y - F^*_y]
      \leq \rho^K\, \mathbb E_y[\Phi_0] + \tfrac{1+\alpha}{\alpha} \Bigl( \textstyle\sum_{k=1}^K \rho^{K-k} \mathbb E_\mathbb Q[\hat\Delta^k] + \tfrac{3M}{4(1-\rho)}\sqrt{(\mathrm{KL} + C)/N} \Bigr).
    \end{equation}
\end{theorem}
The three terms of \cref{eq:probconv} align with the learned scheme.
The initial Lyapunov energy is independent of the networks and decays at the linear rate $\rho$.
Each term in the sum of discounted training risks is the empirical risk $\hat\Delta^k$ of \cref{eq:empiricalrisk} averaged over $\mathbb Q$.
Early iterations are discounted as discussed in \cref{remark:discount}.
The last term on the right-hand side of \cref{eq:probconv} depends on the trained weights only through the parameter drift in $\mathrm{KL}$, which warm-starting in \cref{alg:layerwisetraining} helps to keep small.\looseness=-1

To run the scheme beyond a trained depth $K_0$ fixed before training, we reuse the last trained networks, $\Theta^k = \Theta^{K_0}$ for additional iterations $k > K_0$, so the divergence stays at its value $\mathrm{KL}_{K_0}$.
Splitting the discounted sum in \cref{eq:probconv} at $K_0$ then gives, for any fixed $K \geq K_0$ (see \cref{cor:probconv}),
\begin{equation}
  \mathbb E_\mathbb Q \mathbb E_y[\tilde F^K_y - F^*_y]
  \leq \frac{1+\alpha}{\alpha(1-\rho)} \Bigl( \hat\Delta^{>K_0} + \tfrac{3M}{4}\sqrt{(\mathrm{KL}_{K_0} + C)/N} \Bigr) + \rho^{K-K_0} B_{K_0},
\end{equation}
where $\hat\Delta^{>K_0} := \max_{K_0 < k \leq K} \mathbb E_\mathbb Q[\hat\Delta^k]$ is the largest training risk of the reused networks and $B_{K_0}$ depends only on the first $K_0$ iterations.
The second term decays at the rate $\rho$ in $K$, so for large horizons the first term dominates.
This term matches the neighborhood in case iii.\ of \cref{corollary2}, but the unknown error bound $\Delta^\text{tail}$ is now replaced by quantities computed from the training data.
It stays small as long as the reused networks remain accurate on the subproblems beyond the training horizon $K_0$.\looseness=-1

\section{Dictionary regularized CT image reconstruction}
\label{sec:ctresults}
The motivating application of our learned optimization approach is the reconstruction of computed tomography (CT) images using dictionary-based regularization in low-dose settings \citep{rudzusika2022deep}.
For this particular use-case, the objective function has the form of \cref{eq:Objective} and is given by\footnote{To obtain the form of \cref{eq:Objective}, we rescale the objective by $1/(2\lambda_1)$ which does not change the minimizer.}
\begin{equation}
    \label{eq:objective-ctrecon}
    F(x, z; y) = \| T x - y \|_{\mathcal{Y}}^2 + \lambda_1 \| Dz - x \|_{\mathcal{X}}^2 + \lambda_2 \| z \|_1
\end{equation}%
where the linear operator $T: \mathcal{X} \to \mathcal{Y}$ models how noise-free tomographic data is generated (forward operator), and $D: \mathcal{Z} \to \mathcal{X}$ is the linear synthesis operator for the fixed dictionary that reconstructs the image from the dictionary coefficients, which is implemented as a convolution operator.
We refer to the optimization variable $x$ as the reconstructed image and $z$ as the dictionary coefficients.
The optimization problem is parameterized by the problem data $ y$ given as noisy indirect observations.\looseness=-1

\begin{proposition}
    \label{prop:qg}
    For every bounded set of data $\mathcal{Y}_0$ and every bounded set $\mathcal{S} \subseteq \mathcal{X} \times \mathcal{Z}$, there is a constant $\mu > 0$ such that, for every $y \in \mathcal{Y}_0$, the dictionary-regularized CT objective from \cref{eq:objective-ctrecon} satisfies the quadratic growth condition (\cref{assumption:quadraticgrowth}) at every point of $\mathcal{S}$.\looseness=-1
\end{proposition}

Assuming $p_y$ has bounded support and the learned iterates lie in a bounded set (\cref{sec:learnedopt}), \cref{prop:qg} gives one constant $\mu$ for all instances, as required in \cref{sec:generalization}.
The constant depends on a Hoffman constant that cannot be computed at this scale, so we estimate $\mu$ numerically along the learned iterates
in \cref{fig:qgassumption} (\cref{sec:verifications}).
We now compare the learned scheme to the inexact baseline that solves each subproblem with first-order methods with varying budgets for the number of inner iterations.\looseness=-1

\begin{wrapfigure}{r}{0.5\textwidth}
    \vspace{-2em}
    \centering
    \includegraphics[width=\linewidth]{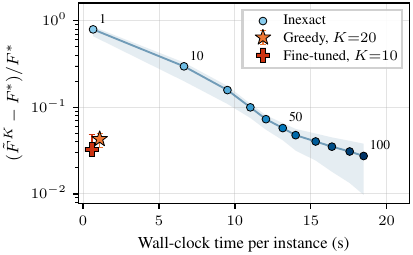}
    \caption{Runtime vs.\ final relative objective gap of learned and baseline methods: medians with 5th--95th percentiles and numbers give the inner iterations of the inexact baseline.}
    \label{fig:pareto}
\end{wrapfigure}
\subsection{Learned vs.\ inexact solvers}%
We evaluate two variants of the learned scheme.
The first uses the greedy training of \cref{sec:greedytr}, which directly targets the per-iteration risk that \cref{theorem:linearbound} bounds.
The second approach applies additional end-to-end fine-tuning to the greedily trained networks.
\Cref{fig:pareto} shows the trade-off between the wall-clock runtime and the final relative objective gap for $20$ outer iterations, traced out as a Pareto front by sweeping the inner-iteration budget of the inexact baseline from $1$ to $100$ inner iterations.
Both learned variants lie well below the frontier traced by the inexact baseline.
Over $5$ independent training runs, the greedy variant reaches a relative objective gap of $0.043 \pm 0.001$ (mean $\pm$ standard deviation) in $1.1$ seconds.
The inexact baseline needs roughly $60$ inner iterations and $14$ seconds to match this gap, a $13\times$ runtime reduction.
The fine-tuned variant lowers the gap further to $0.033 \pm 0.001$ in $0.54$ seconds.
Matching this performance requires roughly $75$ inner iterations and $16$ seconds, leading to an overall $30\times$ reduction.

\subsection{Convergence of alternating scheme}
\begin{figure}
    \centering
    \includegraphics[width=\linewidth]{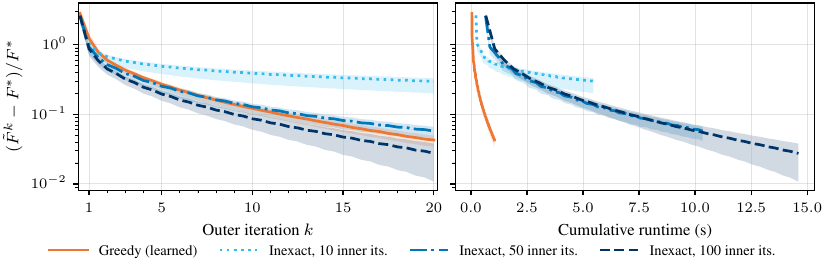}
    \caption{Relative objective gap vs.\ outer iteration $k$ (left) and cumulative runtime (right) over $20$ iterations: medians with 5th--95th percentiles over test instances.}
    \label{fig:convergence}
\end{figure}

In \cref{fig:convergence} we show the per-iteration relative objective gap for the greedy learned scheme and three inexact baselines as a function of the outer iteration $k$ (left) and the cumulative runtime (right).
After the first few outer iterations, all four trajectories decay roughly linearly on the log scale in line with the bound derived in \cref{theorem:linearbound}.
Plotted against wall-clock time, the trajectories of the inexact baselines span an order of magnitude more time, since they are dominated by the inner solver times.\looseness=-1

In the first few outer iterations, the learned scheme is slightly worse in terms of objective value than the inexact baseline with $50$ inner iterations and only catches up around $k \approx 8$.
This is consistent with \cref{remark:discount}: errors in early iterations enter the bound with a discount factor and therefore tolerate larger magnitudes.
The convergence of the learned scheme is therefore governed by the approximation error of the neural network, which we evaluate empirically in the following section.

\subsection{Solving of subproblems}
Here, we analyze how well the greedily-trained networks $\phi$ and $\psi$ solve their respective subproblems from the inexact scheme \cref{scheme2} in isolation for each iteration $k$.
In other words, we are examining the performance independently of error accumulation across iterations.
Given the iterates from \cref{scheme:learned}, we evaluate the function gaps $\xi^k_y$ and $\zeta^k_y$ defined in \cref{eq:inexact-error} by comparing the network predictions to an approximate ground truth solution obtained from a first-order method starting at the current iterate.\looseness=-1
\begin{figure}[b]
    \centering
    \includegraphics[width=\linewidth]{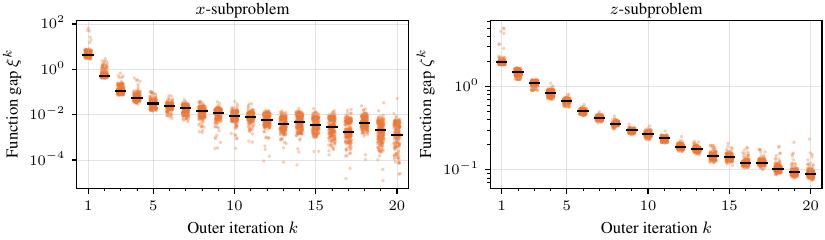}
    \caption{Network approximation error of both subproblems across iterations $k$; dots are test instances and black bars medians.}
    \label{fig:outercomparison}
\end{figure}

The empirical performance of the greedily trained network on the test data is shown in \cref{fig:outercomparison}.
On the left plot, the $x$-subproblem in the first iteration has a larger error which rapidly decreases.
The error for the $z$-subproblem decreases as well on the right plot but the decay is slower and the magnitude of the error is larger in later iterations.
This is unsurprising: the $z$-subproblem in \cref{eq:objective-ctrecon} involves a non-smooth $\ell_1$ term and contains significantly more variables, making it an inherently harder optimization problem.
However, both networks are able to solve the subproblem to a high accuracy and improve the subproblem solution quality with additional outer iterations.\looseness=-1

\subsection{Comparison to other learning approaches}
\label{sec:learnedcomparison}
Additionally, we compare both of our variants to two established learned optimization methods without convergence guarantees.
The first approach uses a \textit{fully-amortized} network to predict the solution.
In other words, the network is trained to directly predict the optimal solution to the joint optimization problem in \cref{eq:Objective} in a single forward pass \citep{amos2023tutorial}.
The second approach uses a \textit{semi-amortized} approach based on algorithm unrolling \citep{monga2021algorithm}.
The network uses the same architecture as our proposed method but we apply a forward pass through all layers combined and train the model end-to-end.
More details about the model implementations can be found in \cref{sec:learneddetails}.\looseness=-1

The fine-tuned greedy variant achieves the lowest test objective in \cref{tab:learnedcomp}.
It reaches $23.84$ at $K = 10$ while greedy alone reaches $24.04$ at $K = 20$, both below the $24.43$ of the best baseline.
The fully-amortized baseline performs substantially worse, confirming that one-shot prediction is the wrong approximation class for this problem.
The end-to-end unrolling approach cannot be extended to more iterations because backpropagation through more unrolled iterates exceeds GPU memory.
In contrast, greedy training optimizes one iteration at a time and so scales freely with $K$.
At its maximum supported depth ($K = 10$), end-to-end unrolling reaches a lower objective than greedy alone ($24.43$ vs.\ $25.82$), but greedy trains three times faster ($11$h$21$m vs.\ $34$h$11$m).
Adding a few steps of end-to-end fine-tuning to the greedy networks closes this gap and produces the lowest objective overall ($23.84$), at roughly half the training cost of the unrolled baseline.

\input{tables/learned_comparison}

\section{Robust matrix completion}
\label{sec:matrixresults}
A second, synthetic, example is the robust matrix completion problem \citep{tao2011recovering}, also known as low-rank matrix completion with sparse corruptions.
Given partial observations of an $n \times m$ matrix $y$ that is the sum of a low-rank and a sparse matrix, as for example in video sequences, we aim to recover both components by minimizing \citep{klopp2017robust}
\begin{equation}
    \label{eq:objectiverobust}
    F(x, z; y) = \lambda_1 \| x \|_* + \frac12 \| P_\Omega(y - x - z)\|^2_F + \lambda_2 \| z  \|_1,
\end{equation}
where the optimization variables $x$ and $z$ are $n \times m$ matrices.
The sampling operator $P_\Omega$ keeps the entries in the observed index set $\Omega$ and sets all other entries to zero.
The nuclear norm favors a low-rank component $x$, the $\ell_1$ norm promotes a sparse $z$, and the masked Frobenius norm couples the two variables.
The regularization parameters $\lambda_1$ and $\lambda_2$ are chosen to minimize the reconstruction error of the minimizer with respect to the ground truth.
We consider $100 \times 100$ matrices of rank $5$ with $40\%$ of the entries observed and $5\%$ of the entries corrupted by sparse noise (\cref{sec:matrixdata}).

\begin{proposition}
    \label{prop:qg2}
    If a minimizer of \cref{eq:objectiverobust} satisfies the standard robust PCA assumptions of strict complementarity and restricted injectivity, then it is the unique minimizer and the robust completion objective satisfies the quadratic growth condition (\cref{assumption:quadraticgrowth}) on every sublevel set.\looseness=-1
\end{proposition}
The constant in this proposition depends on the problem instance.
In \cref{app:qg2-common} we obtain a common constant on an event of instances and combine the probability of its complement with \cref{eq:probconv}.\looseness=-1

\myparagraph{Methods}
The $z$-subproblem has a closed-form solution, since it is separable in the entries of $z$.
The alternating scheme therefore computes the $z$-update exactly, so the $z$-error $\zeta^k_y$ in \cref{eq:inexact-error} vanishes.
The $x$-subproblem has no closed-form solution, and each proximal gradient step requires singular value thresholding (SVT) with a computationally expensive singular value decomposition (SVD).
In our learned method, the network unrolls a single SVT step in the style of LISTA \citep{gregor2010learning} and learns a step size and one threshold per singular value, that is, $101$ parameters per iteration $k$.
We compare against truncated FISTA, a tuned proximal gradient method, and ADMM, whose reformulation gives closed-form updates with a single SVD per iteration (\cref{sec:matrixbaseline}).
\begin{wrapfigure}{r}{0.5\textwidth}
    \centering
    \includegraphics[width=\linewidth]{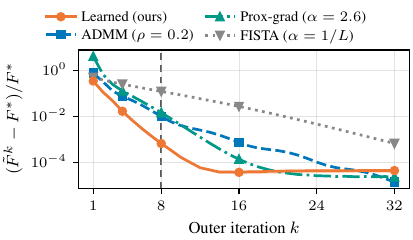}
    \caption{Relative objective gap for robust matrix completion against the outer iteration $k$ with equal computational budget per iteration.
    The learned curve is the mean over $5$ training seeds.
    The dashed vertical line marks the trained depth $k = 8$.}
    \label{fig:rcconvergence}
    \vspace{-0.5em}
\end{wrapfigure}

\myparagraph{Results}
\Cref{fig:rcconvergence} shows the relative optimality gap on $128$ unseen test instances at an equal budget of one SVD per outer iteration.
Up to the trained depth $k = 8$, the learned scheme attains the lowest gap, by a factor of $1.4$ to $15$ over the best baseline (\cref{tab:robustcompletion} in \cref{sec:matrixadditional}).
At $k = 8$, FISTA needs $5\times$ more wall-clock time to reach the objective of the learned scheme (\cref{tab:robustcompletiontiming}).
Beyond the training horizon, we reuse the last trained networks as in \cref{cor:probconv}.
The gap further decreases until $k \approx 14$ and then plateaus, which is consistent with the error term of the late training risk in \cref{cor:probconv}.
The tuned proximal gradient method and ADMM keep decreasing and overtake the learned scheme at $k = 20$ and $k = 29$, respectively.

\section{Discussion \& Conclusion}
\label{sec:conclusion}

\myparagraph{Limitations and future work}
The convergence analysis in \cref{sec:theory} relies on joint convexity with quadratic coupling as in \cref{eq:Objective} and the quadratic growth condition (\cref{assumption:quadraticgrowth}) which excludes bi-convex objectives common in alternating minimization.
A natural next step is to relax these assumptions.
Further, the empirical evaluation covers one imaging modality with a single held-out patient at a fixed noise level and fixed regularization weights $\lambda_1, \lambda_2$, and synthetic robust matrix completion instances, so variability across patients, problem parameters, and real data is not evaluated.
Further scaling the imaging from 2D to 3D CT reconstruction can lead to even larger speedups.\looseness=-1

\cref{thm:pacbayes} assumes that training and evaluation instances are drawn from the same distribution $p_y$, and thus does not cover distribution shift, a standard limitation of learning-to-optimize methods.
Moreover, the bound in \cref{thm:pacbayes} concerns the posterior $\mathbb Q$ at a horizon fixed in advance, while the deployed point estimate is covered only by the per-instance bound of \cref{theorem:linearbound}.
Evaluated with an empirical estimate of $\mu$ and other plug-in constants, the PAC-Bayes bound on the discounted error sum is informative for robust matrix completion but not for CT reconstruction, and the resulting bound on the optimality gap is vacuous for both (\cref{app:pacbayes-ct,app:pacbayes-rc}).\looseness=-1

The largest speedup in the experiments comes from end-to-end fine-tuning, where the training objective no longer coincides with the per-iteration risk in the bound, so the analysis covers these networks but does not motivate their training directly.
Our greedy scheme uses separate weights $\Theta^k$ for every trained iteration, which scales the parameter count and inference-time memory linearly.
Reusing the last trained networks caps this growth at the cost of an error term in the bound that does not vanish.\looseness=-1

\myparagraph{Conclusion}
In this paper, we accelerate alternating minimization schemes of convex objectives with quadratic coupling by replacing each subproblem with an inexact solver that can be implemented by a neural network.
For the inexact scheme, \cref{theorem:linearbound} gives a worst-case bound on the optimality gap, and \cref{corollary2} characterizes three convergence regimes dependent on the decay of the per-iteration error.
Building on this bound, we train a sequence of neural networks greedily by iteration and combine the inexact analysis with a PAC-Bayes argument to obtain a fixed-horizon high-probability guarantee on the expected optimality gap under the posterior of the trained network (\cref{thm:pacbayes}).

For the dictionary-regularized CT reconstruction, the learned scheme traces out a runtime/accuracy frontier strictly below that of the inexact baseline and improves on both fully-amortized and unrolled learned alternatives.
In the synthetic robust matrix completion problem, our method outperforms classical solvers at the trained depth under the same computational budget and remains ahead of them up to twice the trained depth.

\clearpage
\begin{ack}
We would like to thank the anonymous reviewers for their helpful feedback on earlier versions of this paper.
This work was supported by the Wallenberg AI, Autonomous Systems and Software Program (WASP) funded by the Knut and Alice Wallenberg Foundation.
\end{ack}

\bibliography{bibliography}

\newpage

\appendix

\section{Preliminaries}

\subsection{Basic notions and facts from convex analysis}
We review some key results from optimization and convex analysis used in the proofs.
The main reference for this part is \citep{nesterov2018lectures} but some results can also be found in \citep{BauschkeCombettes:2017}.
In the following, let $\mathcal{X}$ be a Hilbert space with inner product $\langle \cdot, \cdot \rangle$ and induced norm $\| \cdot \|$.
Throughout this paper we typically use $\mathcal{X} = \mathbb{R}^n$.\looseness=-1

We next list some basic notions and facts from convex analysis that are used throughout the paper. When stating these, $f \colon \mathcal{X} \to \mathbb{R}$ will be a functional that is, unless stated otherwise, assumed to be differentiable. 
Finally, $x, y \in \mathcal{X}$ are two arbitrary points in the domain of $f$.
\begin{itemize}
    \item Assume $f$ is $L$-smooth, i.e., its gradient $\nabla f$ is Lipschitz continuous with constant $L$. Then we have the following inequality (\emph{descent lemma}):
\begin{equation} \label{eq:descentlemma}
    f(y) \leq f(x) + \bigl\langle \nabla f(x), y-x \bigr\rangle + \frac{L}{2} \|y-x\|^2.
\end{equation}

    \item A differentiable function $f$ is convex if and only if
\begin{equation} \label{eq:conv}
    f(x) \geq f(y) + \bigl\langle \nabla f(y), x-y \bigr\rangle.
\end{equation}
    \item Furthermore, $f$ is strongly convex with modulus $\mu > 0$ if and only if 
    \begin{equation}
        \label{eq:strong_conv} 
        f(x) \geq f(y) + \langle\nabla f(y), x-y\rangle + \frac{\mu}{2} \|y-x\|^2.
    \end{equation}

    Equivalently, $f$ is $\mu$-strongly convex iff  $x \mapsto f(x) - \frac{\mu}{2}\|x\|^2$ is convex.

    \item If $f$ is non-smooth but convex, then the set of subgradients is defined as follows; $g$ is a sub-gradient of a convex function, denoted by $g \in \partial f(y)$, if
    \begin{equation} \label{eq:subdif}
        f(x) \geq f(y) + \langle g, x-y\rangle \quad \text{for any} \quad x.
    \end{equation}
    In other words, the local approximation of $f$ at the point $y$ through $g$ is a global lower bound of $f$.
    If $f$ is Fr\'echet differentiable, then $g = \nabla f(y)$.

    \item A proper lower semi-continuous function is closed. 
    \item If $f$ and $g$ are closed convex functions then  $\partial (f+g) = \partial f + \partial g $.
    
    \item If $W \colon \mathcal{X} \to \mathcal{X}$ is a bounded, self-adjoint, positive definite linear operator, then
    \begin{equation}\label{eq:WNorm}
        \| x \|_W^2 := \langle Wx, x \rangle = \bigl\langle W^{\frac{1}{2}} x, W^{\frac{1}{2}} x \bigr\rangle
    \end{equation}
    defines a norm on $\mathcal{X}$.
    If $W$ is only positive semi-definite, then the above is a semi-norm as positivity is not guaranteed since the elements can be mapped to the nullspace of the operator $W$ that parametrizes the norm.
    \item The arithmetic mean - geometric mean (AM-GM) inequality (also referred to as Young's inequality) states that for any $a, b \in \mathbb{R}$ and any $\lambda > 0$:
    \begin{equation}
    ab \leq \frac\lambda2a^2 + \frac{1}{2\lambda}b^2.
    \end{equation}
    \item A simple but useful fact about norms is the following quadratic expansion of squared norms:
    \begin{equation}
        \frac12 \| a + b \|^2 - \frac12 \| a\|^2 = \frac12 \| b \|^2 + \langle a, b \rangle.
    \end{equation}
    Further, the parallelogram law which is satisfied by norms that are induced by an inner product:
    \begin{equation}
        \| x + y \|^2 + \| x - y \|^2 = 2\|x\|^2 + 2\|y\|^2.
    \end{equation}
\end{itemize}

\subsection{Auxiliary lemmas for the alternating scheme}
\label{sec:background}
In this section, we provide some background on the convergence of alternating schemes and provide some key technical lemmas.

The key step in proving the convergence of the alternating scheme relies on a generalized notion of strong convexity that we state below, 
which allows us to couple the updates across the two variable blocks through the quadratic term in \cref{eq:Objective}.

\begin{lemma}[Generalized strong convexity]
    \label{lemma2}
    Assume $x \mapsto f(x) - \frac{1}{2}\|\matr x\|^2$ for an arbitrary bounded linear operator $\matr$ and norm $\| \cdot \|$ on 
    the Hilbert space of the co-domain is a convex function and $g \in \partial f(y)$.
    Then, for any $x$
\begin{equation} \label{eq:gen_strong_conv}
    f(y) \leq f(x) + \langle g, y - x \rangle - \frac{1}{2}\|\matr (y-x)\|^2.   
\end{equation}
If $y$ is additionally a minimizer of $f$, then 
\begin{equation} \label{eq:gen_strong_conv_min}
    f(x) \geq f(y) + \frac{1}{2}\|\matr (y-x)\|^2.
\end{equation}
\end{lemma}
The proof mimics the strongly convex case using instead the semi-norm induced by the linear operator $A$ to weight the norm.
\begin{proof}
Assume that $g \in \partial f(y)$ is a subgradient of $f$, then by applying subdifferential calculus $g - \matr^T\matr y \in \partial (f(y)  - \frac{1}{2}\|\matr y\|^2) $ and by using the subgradient property \cref{eq:subdif}, we obtain
\begin{equation}
\begin{split}
    f(x)  - \frac{1}{2}\|\matr x\|^2 &\geq f(y)  - \frac{1}{2}\|\matr y\|^2+ \langle g - \matr^T\matr y, x-y \rangle \\
    & = f(y)  - \frac{1}{2}\|\matr y\|^2 + \langle g, x-y \rangle -  \langle \matr y, \matr x \rangle + \|\matr y\|^2 \\
    & = f(y)  + \frac{1}{2}\|\matr y\|^2 -  \langle \matr y, \matr x \rangle + \langle g, x-y \rangle.
\end{split}
\end{equation}
By rewriting and exploiting the linearity of the inner product, we obtain the following inequality
\begin{equation}
\begin{split}
    f(x)  &\geq f(y)  + \frac{1}{2}\|\matr x\|^2 + \frac{1}{2}\|\matr y\|^2 -  \langle \matr y, \matr x \rangle  + \langle g, x-y \rangle  \\
    & = f(y) + \frac{1}{2}\|\matr (x- y)\|^2+  \langle g, x-y \rangle .
\end{split}
\end{equation}
If $y$ is the optimum, then $g = 0 \in \partial f(y) $. If $f$ is differentiable $\partial f(y) = \{ \nabla f(y) \}$.
\end{proof}

We can immediately observe that the objective function \cref{eq:Objective} satisfies the assumption of \cref{lemma2} in the following remark. 
\begin{remark}
\label{remark1}
    Let $F$ be a function of the form given in \cref{eq:Objective}.
    Then the functions $x \mapsto F(x, z) - \frac{1}{2}\|A x \|^2$ and $z \mapsto F(x, z) - \frac{1}{2}\|B z \|^2$ are convex.
\end{remark}
\begin{proof}[Proof sketch]
    This can be easily verified as the following is convex in $x$ for a fixed $z$:
\begin{equation}
\begin{split}
    F(x,z) - \frac{1}{2}\|A x\|^2
    = f(x)  + \frac{1}{2}\|A x + B z \|^2 - \frac{1}{2}\|A x\|^2 + g(z) \\
    = f(x)  + \langle A x, B z \rangle + \frac{1}{2}\|B z \|^2 + g(z),
\end{split}
\end{equation}
since $f$ is assumed to be convex and the inner-product term is linear in $x$.
The convexity in $z$ follows by the same argument.
\end{proof}

In other words, the subproblems are generalized strongly convex.
Now, we have all ingredients to prove the following lemma that allows us to derive the convergence of the alternating scheme where the subproblems are solved exactly:
\begin{lemma}
    \label{lemma0}
    Let $F$ be a convex objective that is block separable with quadratic coupling term as in \cref{eq:Objective}.
    Then,  it holds for the iterates obtained from the alternating scheme \cref{scheme1} and any $(x, z) \in \mathcal{X} \times \mathcal{Z}$ that
    \begin{equation}
        \label{eq:lemma0eq}
        F(x^{k+1}, z^{k+1} ) + \frac{1}{2}\|B(z - z^{k+1})\|^2 \leq F(x,z) + \frac{1}{2}\|B(z - {z}^{k})\|^2
    \end{equation}
\end{lemma}
\begin{proof}
    Using the observation from \cref{remark1} and applying \cref{lemma2} yields the following two inequalities 
    for the updates of the alternating scheme \cref{scheme1}. 
    Here, we are using the fact that the next iterates themselves are computed as the solution to an optimization problem:
\begin{subequations}
\label{eq:two_ineq}
\begin{align}
    F(x^{k+1}, {z}^k ) + \frac{1}{2}\|A(x - x^{k+1})\|^2 & \leq F(x, {z}^k ), \label{eq:alt_ineq1}\\
    F(x^{k+1}, z^{k+1} ) + \frac{1}{2}\|B(z - z^{k+1})\|^2 & \leq F(x^{k+1}, z ),\label{eq:alt_ineq2} 
\end{align}
\end{subequations}
which hold for any $x \in \mathcal{X}$ and $z \in \mathcal{Z}$.
Summing these inequalities and using an explicit form of the objective function \cref{eq:Objective} results in
\begin{equation}
\label{eq:summedineq}
\begin{split} 
    F(x^{k+1}, z^{k+1} ) + \frac{1}{2}\|B(z - z^{k+1})\|^2 &\leq F(x,z) + \frac{1}{2}\|B(z - {z}^{k})\|^2 + \mathcal{R}^{k+1}(x, z)
\end{split}
\end{equation}
where the residual can be simplified and bounded in the following way 
\begin{equation}
\label{eq:residualnonneg}
\begin{split}
    \mathcal{R}^{k+1}(x, z) &= \frac{1}{2} \|Ax + B{z}^k\|^2  +\frac{1}{2} \|Ax^{k+1} + B z\|^2 - \frac{1}{2}\|Ax^{k+1} +  B{z}^k \|^2 \\
    &\qquad - \frac{1}{2}\|A x + Bz\|^2- \frac{1}{2}\|A(x - x^{k+1})\|^2 - \frac{1}{2}\|B(z - {z}^{k})\|^2 \\
    &= -\frac{1}{2} \|Ax^{k+1} + B {z}^{k} - A x - Bz\|^2 \leq 0.
\end{split} 
\end{equation}
Plugging this back into \cref{eq:summedineq} finalizes the proof.
\end{proof}

In \cref{sec:theory} we extend the alternating scheme to inexact solves, making use of \cref{lemma2,remark1}.
In particular, we extend the inequality in \cref{eq:lemma0eq} to account for the inexactness error introduced by solving each problem inexactly.

\section{Analysis of the inexact scheme}
\label{sec:proofs}
This section contains the proofs of the results for the inexact scheme stated in \cref{sec:theory}.
We provide the proofs for \cref{lemma0p}, \cref{theorem:linearbound}, and \cref{corollary2} in this appendix.

\subsection{Proof of \cref{lemma0p}}
The proof closely follows the proof of the exact case in \cref{lemma0}, adds error terms from the inexact solvers, and bridges
the exact $z$-minimizer on the left-hand side to the inexact iterate $\tilde{z}^{k+1}$.
The key difference of \cref{eq:lyapanov} compared to the exact case is the additional combined error term and the inflated coefficients on the norm terms.

\input{proofs/lemma_1}

\subsection{Proof of \cref{theorem:linearbound}}
The proof of the linear convergence result relies on the inexact descent inequality shown in \cref{lemma0p} and the quadratic growth assumption formulated in \cref{assumption:quadraticgrowth}.
The final key step is determining the bounds for $\alpha$ such that the rate $\rho < 1$ to obtain convergence of the scheme.

\input{proofs/theorem_1}

\subsection{Proof of \cref{corollary2}}
Since the first term in \cref{eq:theorem1} decays geometrically and the weights $\rho^{K-k}$ discount early errors exponentially, the behavior of the bound for large $K$ is dictated by how fast the late errors $\Delta_k$ decay, which leads to the three regimes listed in \cref{corollary2}.

\input{proofs/corollary_1}

\section{Analysis of the learned scheme}
\label{sec:proofs-learned}
This section contains the proof of \cref{thm:pacbayes} and provides the formal statement for the weight reuse that is mentioned in the main text (\cref{app:reuse}).
Then, we describe how we evaluate the bound numerically (\cref{app:pacbayes-eval}).

\subsection{Proof of \cref{thm:pacbayes}}
\input{proofs/theorem_2}

The prior variance $\varsigma^2$ enters the prior, and the constants $\alpha$, $\mu$ (through the contraction rate $\rho$), $M$, and the horizon $K$ enter the loss $L_y$, so together with $\beta$ and $J$ they must be fixed before training.

If $\varsigma$ and $\alpha$ are selected from a finite grid of $G$ values, applying the argument above at each grid point with confidence $\beta/G$ and taking a union bound shows that \cref{eq:probconv} holds simultaneously for all grid points where the constant $C$ needs to be replaced by $C + \log G$.
The inequality \cref{eq:pacbayesraw} holds simultaneously for all posteriors $\mathbb Q$.
Choosing among several trained models afterwards, all with the same initialization $\Theta^0$ and the same grid, is therefore covered by the same event and requires no further union bound.

\subsection{Weight reuse beyond the trained depth}
\label{app:reuse}
This section contains the bound for weight reuse from \cref{sec:generalization}.
The scheme runs for $K$ outer iterations beyond a trained depth $K_0$ which is fixed before training, by reusing the last trained networks, $\Theta^k = \Theta^{K_0}$ for $k > K_0$.
\Cref{cor:probconv} splits the bound of \cref{thm:pacbayes} at an index $m$ into an asymptotic term, which depends on both the training risk after iteration $m$ and the divergence $\mathrm{KL}_{K_0}$, and a transient term, which decays linearly with the rate $\rho^{K-m}$.
The split holds simultaneously for all $1 \leq m \leq K$, and the main text uses $m = K_0$.

\begin{corollary}
    \label{cor:probconv}
    Let the assumptions of \cref{thm:pacbayes} hold.
    We fix a horizon $K_0$ before training such that $K_0 \leq K$ and reuse the last trained networks for $k > K_0$.
    
    For a split index $1 \leq m \leq K$, define the late training risk $\hat\Delta^{> m} := \max_{m < k \leq K}\mathbb E_\mathbb Q[\hat\Delta^k]$, with $\hat\Delta^{> m} := 0$ if $m = K$, and the transient constant of the early iterations
    \begin{equation*}
        B_m := \rho^{m}\, \mathbb E_y[\Phi_0] + \frac{1+\alpha}{\alpha} \sum_{k=1}^{m} \rho^{m-k}\, \mathbb E_\mathbb Q[\hat\Delta^k].
    \end{equation*}
    Then, for any $\beta \in (0,1)$, with probability at least $1-\beta$, it holds for all $1 \leq m \leq K$ that
    \begin{equation}
        \label{eq:probconvfinite}
        \mathbb E_\mathbb Q \mathbb E_y[\tilde F^K_y - F^*_y]
        \leq \frac{1+\alpha}{\alpha(1-\rho)} \Bigl( \hat\Delta^{>m} + \tfrac{3M}{4} \sqrt{(\mathrm{KL}_{K_0} + C)/N} \Bigr) + \rho^{K-m} B_m,
    \end{equation}
    where $\mathrm{KL}_{K_0}$ is the divergence of the random walks on $\Theta^{1:K_0}$.
\end{corollary}
\input{proofs/corollary_2}

\Cref{cor:probconv} holds at a single horizon $K$ fixed in advance.
For finitely many horizons $K \leq K_\mathrm{max}$ fixed in advance a simple union bound argument already gives a simultaneous statement: applying \cref{cor:probconv} at each horizon with confidence $\beta / K_\mathrm{max}$ replaces $C$ by $C + \log K_\mathrm{max}$.
Since $C$ enters only through $\sqrt{\mathrm{KL} + C}$, this penalty is small whenever $\mathrm{KL}$ dominates.

\subsection{Numerical evaluation of the bound}
\label{app:pacbayes-eval}
This section describes how we evaluate \cref{thm:pacbayes} numerically for both problems.
The problem-specific settings and results are given in \cref{app:pacbayes-ct} for CT reconstruction and in \cref{app:pacbayes-rc} for robust matrix completion.

We evaluate the inequality \cref{eq:pacbayes} from the proof of \cref{thm:pacbayes}, which bounds the discounted posterior risk $\sum_k \rho^{K-k} \mathbb E_\mathbb Q[\mathbb E_y[\Delta^k_y]]$ directly.
The full theorem statement \cref{eq:probconv} multiplies this quantity by $(1+\alpha)/\alpha$ and adds the decaying initial term $\rho^K \mathbb E_y[\Phi_0]$.
Since $\alpha$ enters \cref{eq:pacbayes} only through $\rho$, the selected $\alpha$ lies at the lower end of the grid where $(1+\alpha)/\alpha$ is large.
We therefore report the evaluation as a bound on the discounted error sum itself.

\myparagraph{Plug-in constants}
Two constants of the theorem are not available in certified form, and we replace them by plug-in estimates we obtain from data.
\begin{itemize}
    \item \emph{The bound $M$.}
    The theorem requires $\Delta^k_y \leq M$ for every instance and every parameter value.
    The output projection guarantees that such an $M$ exists, but the resulting a priori value is several orders of magnitude larger than the errors that occur and makes the bound vacuous on both problems.
    We therefore use $\hat M(\varsigma)$, the largest $\Delta^k_y$ observed over the posterior samples at the variance $\varsigma$.
    We report a running maximum over increasing $\varsigma$, so that $\hat M$ is nondecreasing in $\varsigma$, and bound it from below by the largest $\Delta^k_y$ at the trained weights, computed over the full training set for CT and over the training and test sets for robust matrix completion.
    \item \emph{The constant $\mu$.}
    The contraction rate $\rho$ depends on the quadratic growth constant.
    The theorem requires one constant $\mu$ for all instances and for every parameter value in the support of $\mathbb Q$.
    We use a common empirical constant $\hat\mu$, the smallest quadratic growth ratio observed along the iterates.
    For CT, we compute it along the deployed iterates on the test set (\cref{sec:verifications}).
    For robust matrix completion, we additionally include the iterates under $\mathbb Q$ and the iterates of the classical solvers (\cref{app:pacbayes-rc}).
\end{itemize}
The resulting numbers are a numerical evaluation of the bound, not a certificate.

\myparagraph{Monte Carlo estimates and grid}
The posterior risks $\mathbb E_\mathbb Q[\hat\Delta^k]$ are Monte Carlo estimates over sampled trajectories $\Theta^{1:K} \sim \mathbb Q$.
We do not include the Monte Carlo error in the bound.
We select $\alpha$ and $\varsigma$ from a grid of $G$ fixed values, which by a union bound adds $\log G$ to the constant $C$ in \cref{thm:pacbayes}.
We use $J = 40$ for the grid over the free parameter $\lambda$, which covers every posterior width we evaluate, and $\beta = 0.05$ for both problems.
We compare the numerical evaluation of the PAC-Bayes bound to the trivial bound $\hat M (1 - \rho^K)/(1-\rho)$, which follows from $\Delta^k_y \leq \hat M$ alone.

\section{Dictionary regularized CT reconstruction}
This section contains the proof of \cref{prop:qg}, the data generation, the baseline and learned methods, and additional results for the CT reconstruction experiments, including the evaluation of the PAC-Bayes bound.

\subsection{Proof of \cref{prop:qg}}
\input{proofs/proposition_1}

\subsection{Problem parametrization}
\label{subsec:problem}
\myparagraph{Imaging data}
We are using the computed tomography (CT) imaging data of the human abdomen from the 2016 Low Dose CT Grand Challenge as ground truth \citep{mccollough2021low}.
The dataset contains imaging data from a total of 10 patients out of which we are using the data from 9 patients for training and validation while the final patient is reserved for testing.

\myparagraph{Data simulation}
In this paper, we are only focusing on 2D image reconstruction by splitting the 3D volumetric data into axial slices.
We are using the Operator Discretisation Library (ODL) \citep{adler2017operator} with the astra-toolbox CT backend \citep{van2015astra, van2016fast} to simulate 2D projection data from the extracted $512 \times 512$ image slices.
The acquisition geometry is modeled as a fan beam geometry with $1\thinspace000$ uniformly partitioned angles and $1\thinspace000$ detector elements.
The noise free tomographic data modeled by the forward operator is given by $y = Tx$, where $T$ is the linear Radon transform operator also referred to as Ray Transform.\footnote{Note that $T$ is implemented as a matrix-free operator. It is not possible to form the matrix $T$ directly since the size of the matrix is prohibitively large containing over $262$ billion elements. The matrix-free version directly computes the matrix-vector product by solving line integrals.}

We model the observed noisy data $y_{\text{noisy}} \approx y$ as low-dosage clinical CT data.
Therefore, we re-scale the images given in Hounsfield scale by $\mu_0 / 1\thinspace000$ for the water attenuation constant $\mu_0 = 0.0192$mm$^{-1}$.
For the single energy X-Ray noise model we assume a Poisson distribution with mean corresponding to the number of photons absorbed by the detector. This number of absorbed photons is modeled by the Beer-Lambert law as
\begin{equation}
    N_d = N_0 \cdot \exp{(-Tx)}.
\end{equation}
The noisy observation is given as a linearization of the data:
\begin{equation}
    y_{\text{noisy}} = - \ln \left( \frac{1}{N_0} \; \text{Poisson}(N_d) \right).
\end{equation}
To avoid taking the logarithm of zero, the photon counts are clipped from below at $1$, so that $y_{\text{noisy}} \leq \ln N_0$.
Since Poisson counts have no upper limit, the data is not bounded from below, although large deviations have negligible probability.
In our experiments we choose $N_0 = 50\thinspace000$ to model the clinical low-dose setting.

From this simulation, we obtain pairs of noisy observations and ground truth images $\{ (x, y_{\text{noisy}})\}_{i=1}^N$ with a total of $2\thinspace147$ training samples, 21 validation samples, and 210 samples from a different individual used for testing purposes.
Training the data-driven optimizer, however, is done in an unsupervised fashion where only the observations $y_\text{noisy}$ are required that parameterize the objective.
Further, note that our primary goal is not to recover the ground truth $x$ as well as possible but rather solve the optimization problem \cref{eq:objective-ctrecon} as fast as possible.
How to design a good regularizer such that the solution to the optimization problem is close to the true solution $x$ is out of scope for this paper.\looseness=-1

\myparagraph{Dictionary regularization}
In our experiments, the optimization objective is given by the data fidelity term, a term corresponding to dictionary representation of the signal as well as a sparsifying $\ell_1$ term on the dictionary coefficients.
The dictionary consists of $128$ dictionary elements of size $16 \times 16$ and are pre-trained using the method described by \citet{rudzusika2022deep}.
The synthesis operator $D$ is a linear operator reconstructing an image from the dictionary coefficients using the learned dictionary elements.
As regularization parameters, we choose $\lambda_1=50$ and $\lambda_2=0.0016$ to minimize the validation error in terms of image reconstruction performance.
In order to reduce the number of dictionary coefficients for the reconstruction, we apply a stride $2$ to the elements of the dictionary.\looseness=-1

\myparagraph{Optimization variables}
The total number of optimization variables in the CT reconstruction objective from \cref{eq:objective-ctrecon} is given by the sum of image variables (for the $x$-subproblem) which is of dimension $512 \times 512 = 262\,144$ and the dictionary coefficients (for the $z$-subproblem) which is $128 \times 256 \times 256 = 8\,388\,608$.
Thus, we have a total of $8\,650\,752$ variables.

\subsection{Baseline methods}
\label{sec:baseline}
Here, we discuss how all methods apply the pre- and post-processing.
Further, we specify the details for the solving procedure for the alternating and inexact scheme and how the iterates are numerically computed.
Further, we provide details how the ground truth solution $(x^*, z^*)$ is approximated.

\myparagraph{Pre- and post-processing}
In order to apply the dictionary regularized variational formulation, we are removing the low frequency components from the data with a low-pass filter in the pre-processing step.
Therefore, in a first step we apply the filtered back projection algorithm (FBP) with a Hanning filter and relative frequency cut-off set to $0.1$.
This step is necessary to allow good regularization performance via the dictionary elements that only model the high frequency signal components \citep{rudzusika2022deep}.
To get the final reconstructed image, we add the previously extracted low frequency structures of the image to the result obtained from the learned alternating minimization algorithm.\looseness=-1

In order to initialize our alternating minimization scheme, we are using the analytical filtered back projection (FBP) method with Hanning filter and relative frequency cut-off $0.75$ to obtain an initial guess $x^0$ that is provided to the algorithm.

The baseline methods utilize the Lipschitz constant of the data fidelity term and the quadratic coupling between the variable blocks to set the step-size.
To compute the Lipschitz constant, we use the power method with $10$ iterations which estimates the largest singular values of the linear operators $T$ (ray transform) and $D$ (dictionary synthesis operator) as part of the pre-processing.

\myparagraph{Solving the subproblems}
In each iteration of the alternating minimization scheme, the two subproblems shown in~\cref{scheme1} have to be solved to compute the next iterates.
For both subproblems the solution obtained in the previous outer iteration can be used to warm-start the solver leading to a significant reduction in solving time for the iterative schemes.

The first sub-problem in \cref{scheme1}, which corresponds to the update of the $x$-variable, is both convex and smooth.
Therefore, we solve the problem using the first-order gradient method with Nesterov momentum for the baseline method.
Since this subproblem involves the forward operator in the CT reconstruction problem, it is still computationally demanding.

The second subproblem is non-smooth due to the $\ell_1$ norm used to obtain sparse dictionary coefficients.
Since the proximal operator is available as a closed form solution by the soft shrinkage operator \citep{beck2017first} we apply the forward-backward splitting method. 
Similar to the first subproblem, we use the corresponding accelerated version of the proximal gradient method (FISTA) \citep{beck2009fast}.\looseness=-1

In practice, we stop the iterative scheme to solve the inner subproblems when the relative change in objective value is less than a predefined threshold (in practice we choose $10^{-5}$).
The maximum number of inner iterations is variable giving rise to the different baseline methods.
We run the outer scheme using a fixed number of iterations ($K=20$) and train a neural network to approximate the solution to the subproblem of each iteration as outlined in \cref{sec:learnedopt} and as we discuss in detail in the next section.

\myparagraph{Ground truth solution}
To evaluate the optimality gap, we approximate a minimizer $(x^*, z^*)$ by running the baseline alternating scheme for up to $200$ outer iterations with up to $50$ inner iterations for each subproblem, stopping each loop once the relative decrease of the objective falls below $10^{-5}$.
Since the subproblems are solved only approximately, the reference value $F(x^*, z^*)$ upper-bounds the optimal value $F^*$, so the reported relative gaps underestimate the true gaps.
All methods are evaluated against the same reference ensuring a fair comparison.
Note that the problem is not strongly convex so the found minimizer may be non-unique.

\subsection{Learned methods}
\label{sec:learneddetails}
In each outer iteration $k$ we parameterize the learned solvers $\phi_{\theta^k}$ and $\psi_{\vartheta^k}$ as small convolutional networks tailored to each subproblem.
Both networks use the same underlying backbone consisting of $3 \times 3$ convolutions with $32$ hidden channels which are interleaved with GroupNorm layers and ReLU activations.
The two approaches only differ in terms of input and output projection.
Further, the $z$-network has additional learned soft-thresholding to incorporate the non-smooth $\ell_1$ penalty into the architecture.
We do not share weights across outer iterations, so each of the $K=20$ iterations has its own pair of networks.

\emph{The $x$-networks}:
The $x$-network takes as input the current iterate $\tilde x^k$ stacked with the two terms of the partial gradient $\nabla_x F(\tilde x^k, \tilde z^k)$ separated into the data-fidelity contribution $2 T^T(T \tilde x^k - y)$ and the coupling contribution $2\lambda_1(\tilde x^k - D \tilde z^k)$, giving $3$ input channels.
The network produces a scalar correction $d_x$ of the same spatial dimension as the image, and the next iterate is formed by augmenting a single preconditioned gradient step with this correction
\begin{equation}
    \tilde x^{k+1} = \tilde x^k - \frac{1}{L_x}\bigl(\nabla_x F(\tilde x^k, \tilde z^k) + d_x\bigr),
\end{equation}
where $L_x$ is the Lipschitz constant of $\nabla_x F(\,\cdot\,, \tilde z^k)$ estimated by power iteration during pre-processing.
The backbone uses $4$ convolutional layers with $32$ hidden channels which results in $29\,127$ trainable parameters per outer iteration.

\emph{The $z$-networks}:
The $z$-network is inspired by the LISTA architecture \citep{gregor2010learning} which unrolls steps of a proximal gradient method into a learned shrinkage map.
Unlike the original LISTA formulation, the network is trained in a fully unsupervised fashion against the per-iteration subproblem objective and we use only a single shrinkage step.
The network takes as input the current sparse code $\tilde z^k \in \mathbb{R}^{128 \times 256 \times 256}$ with $128$ dictionary atoms living on a coarser spatial grid than the image due to the stride-$2$ dictionary $D$, stacked with the smooth-part gradient $2\lambda_1 D^T(D \tilde z^k - \tilde x^{k+1})$, giving $256$ input channels.
The network produces a correction $d_z$, and the next iterate is obtained from an additional proximal-gradient step augmented by $d_z$ followed by atom-wise soft-thresholding
\begin{equation}
      u = \tilde z^k - \frac{1}{L_z}\bigl(g_z + d_z\bigr),\qquad                    
      \tilde z^{k+1} = \text{sign}(u)\odot \max\bigl(|u| - \theta,\, 0\bigr),
\end{equation}
where $L_z = 2\lambda_1 \|D\|^2$ and $\theta \in \mathbb{R}^{128}_{>0}$ is a learned per-atom threshold parameterized as $\theta = \mathrm{softplus}(\rho)$ to enforce positivity.
We zero-initialize the final convolution and initialize $\rho$ such that $\theta = \lambda_2 / L_z$ at the start of training, so the network reduces to one classical ISTA step at initialization.
The $z$-network has $139\,328$ trainable parameters per outer iteration.

\emph{Output projection}:
To keep the learned iterates bounded for every parameter value of the posterior as discussed in \cref{sec:learnedopt}, the output of both networks is projected onto an a priori box,
\begin{equation*}
    \tilde x^{k+1} \leftarrow \Pi_{[-0.1,\,0.1]}(\tilde x^{k+1}) \qquad \tilde z^{k+1} \leftarrow \Pi_{[-4,\,4]}(\tilde z^{k+1})
\end{equation*}
applied elementwise.
The box for $x$ acts on the image after the low-frequency component has been removed (\cref{sec:baseline}) and not on the physical attenuation values.
Since the projection is part of the learned map, it is also applied to every parameter sample drawn from the posterior in \cref{app:pacbayes-ct}.
On the parameters trained from the greedy scheme, the projection is inactive across all test instances.

\myparagraph{Training settings}
At each outer iteration $k$ we train the $x$- and $z$-networks sequentially: we first fit $\phi_{\theta^k}$ against the per-iteration $x$-subproblem objective with $\theta^{k-1}$ as a warm start, and then, with $\theta^k$ held fixed, fit $\psi_{\vartheta^k}$ against the per-iteration $z$-subproblem objective with $\vartheta^{k-1}$ as a warm start (see \cref{alg:layerwisetraining}).
Initializing from the previous iteration (line~3 of \cref{alg:layerwisetraining}) lets us reduce the training budget at later iterations to a fine-tuning regime.
The number of training epochs for each network decays with the outer iteration index $k$ according to the schedule
\begin{equation*}                                                                                                                                                       
      (\text{E}_x^k, \text{E}_z^k) = 
      \begin{cases}
          (10, 50) & \text{for } k = 1, \\
          (5,\, 25)  & \text{for } 2 \leq k \leq 6, \\
          (2,\, 10)  & \text{for } k \geq 7,
      \end{cases}                                
  \end{equation*}
where $\text{E}_x^k$ and $\text{E}_z^k$ are the number of epochs used to train the $x$- and $z$-network at iteration $k$, respectively.
We use a batch size of $4$ for the $x$-subproblem and $8$ for the $z$-subproblem.
As the outer scheme converges, the data distribution presented to subsequent networks stabilizes, making the warm start increasingly effective.

Additional hyper-parameters and training settings are shown in \cref{tab:learneddetails}.

\myparagraph{Learned baselines}
We consider two learned baselines that target the canonical objective \cref{eq:objective-ctrecon} but differ in how the predictor is parameterized.
The \textit{fully-amortized baseline} maps the warm-started variables to a prediction of $(x, z)$ in a single forward pass.
The \textit{semi-amortized baseline} retains the iterative structure of the proposed scheme but trains all per-iteration networks jointly via algorithm unrolling \citep{monga2021algorithm}.
The hyperparameters of the two methods are summarized in \cref{tab:learneddetails}.

\textit{Fully-amortized baseline}:
The fully-amortized network predicts $(x, z)$ in a single forward pass from the initial input variables.
We use a Y-shaped architecture that consists of a shared image-domain residual trunk followed by two thin heads, one for each variable.
The trunk consists of an input projection and $6$ residual blocks, each with two $3 \times 3$ convolutions of $96$ filters interleaved with GroupNorm and ReLU activations and a skip connection across the block.
The $x$-head is a $2$-layer convolutional network that produces an additive correction to the input.
The $z$-head is a $3$-layer network that consumes the trunk features projected into atom space concatenated with $\tilde z^0$, $D^T \tilde x$ and $D^T D \tilde{z}^0$, produces a LISTA-style correction, and applies the same atom-wise soft-thresholding step as in the proposed $z$-network.

\textit{Semi-amortized baseline}: 
The semi-amortized baseline uses the same per-iteration architecture as the proposed greedy method but differs in the training strategy.
Rather than training each pair of networks in isolation against the per-iteration subproblem objective, all networks are trained jointly end-to-end with a single loss  on the final iterate $(\tilde x^K, \tilde z^K)$ propagated through the unrolled scheme.
We use $K = 10$ outer iterations for this baseline.

All baselines use the same validation frequency and early stopping based on the validation performance to avoid overfitting.
\input{tables/hyperparameters}

\subsection{Computational resources}
\label{subsec:resources}
All reported experiments are executed on a single NVIDIA TITAN Xp GPU with $12$~GB of memory on an in-house server with $12$ CPU cores and $62$~GB of RAM.

\myparagraph{Training compute}
A complete layer-wise training run of the learned alternating model takes approximately $17$ hours of wall-clock time on the hardware above.
Additional end-to-end fine-tuning takes around $6$ more hours for the model.
The fully-amortized baseline takes about $17$ hours, and the semi-amortized baseline with $K=10$ about $34$ hours.
We repeated full model training $5$ times for the paper.
The training times are also reported in \cref{tab:learnedcomp}.

\myparagraph{Inference compute}
For all wall-clock timings we repeat each evaluation $3$ times after a warm-up pass and report the mean.
On a single test instance, the exact baseline with $100$ inner FISTA iterations per outer step requires about $18$~s, while the learned model solves the same problem in about $1.1$~s, a speedup of roughly $16\times$.
The fully and semi-amortized baselines as well as the end-to-end fine-tuned model run in $0.3$--$0.6$~s per instance.
For learned methods, every test run additionally solves each inner subproblem to high accuracy with $100$ FISTA iterations in order to report the per-subproblem suboptimality used in \cref{fig:outercomparison}.
This dominates the inference cost and brings a full evaluation pass over the $210$ test instances to roughly one day on the same hardware.
The high-accuracy reference solution $(x^*, z^*)$ used as the denominator of the relative objective gap was computed once and required an additional day of compute.

\myparagraph{PAC-Bayes evaluation}
The posterior Monte Carlo sampling in \cref{app:pacbayes-ct} runs the learned scheme with $K = 20$ for $15$ values of $\varsigma$ and $5$ trajectories each on $10$ training slices, including the reference solves for the per-iteration errors, and takes $4.4$ hours.

\myparagraph{Previous experiments}
The full research project required several times the compute of the reported runs.
During development we trained roughly $20$ additional models that are not used in the paper, mostly with reduced numbers of outer iterations between $5$ and $10$ or with shorter training schedules.

\subsection{Licenses and assets}                                                       
\label{sec:licenses}
\Cref{tab:licenses} summarizes the third-party datasets and software used in this work.
All assets are used in accordance with their respective licenses and within the terms of use.

\begin{table}[h]
    \centering                
    \caption{Third-party assets used in this work.}                                      
    \label{tab:licenses}
    \small
    \resizebox{\textwidth}{!}{%
    \begin{tabular}{llll}
      \toprule                                                                           
      Asset & Type & License & Reference \\ 
      \midrule
      Low Dose CT Grand Challenge & Dataset & CC BY-NC 4.0 / & \citep{mccollough2021low} \\
                                            & & NIH Controlled Data Access Policy & \\
      ODL & Software & MPL-2.0 & \citep{adler2017operator} \\ 
      ASTRA Toolbox & Software & GPL-3.0 & \citep{van2015astra,van2016fast} \\       
      PyTorch & Software & BSD-3-Clause & \citep{paszke2019pytorch} \\           
      \bottomrule 
    \end{tabular}
    }
\end{table}

\subsection{Verification of assumptions}
\label{sec:verifications}
Here, we empirically investigate the ability of the trained neural network to fulfill the assumptions based on which we prove the convergence of the inexact scheme.

\begin{figure}[htbp]
    \centering
    \includegraphics{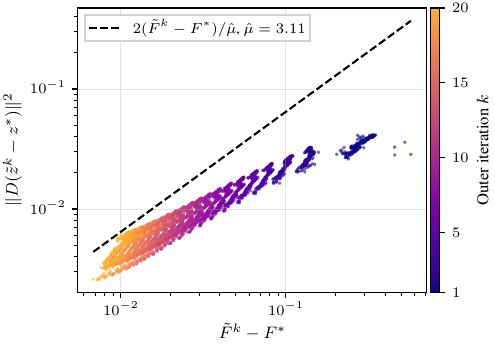}
    \caption{Numerical verification of the assumption of the quadratic growth in a neighborhood of the optimal set as stated in \cref{assumption:quadraticgrowth}.}
    \label{fig:qgassumption}
\end{figure}
\myparagraph{Quadratic growth assumption}
To validate \cref{assumption:quadraticgrowth} we collect the iterates $\{(\tilde x^k, \tilde z^k)\}_{k=1}^{K}$ produced by the learned alternating scheme on each of the $210$ test instances, together with the estimated optimum $(x^*, z^*)$. 
The reference value $F(x^*, z^*)$ and the logged function values $F(\tilde x^k, \tilde z^k)$ are evaluated on the canonical objective \cref{eq:Objective}, that is, after rescaling the objective by $1/(2\lambda_1)$ as discussed in the footnote of \cref{eq:objective-ctrecon}.

We discard the initial iterate $k = 0$ since it precedes any optimization step and is uninformative for the assumption.
For each remaining iterate we plot $\| D(\tilde z^k - z^*)\|^2$ against $F(\tilde x^k, \tilde z^k) - F(x^*, z^*)$ in \cref{fig:qgassumption}, where the color encodes the iteration index $k$.
The empirical quadratic growth constant
\begin{equation}
\hat\mu \;=\; \min_{k,y_i} \; \frac{2\,\bigl(F(\tilde x^k_{y_i}, \tilde z^k_{y_i}; y_i) - F(x^*_{y_i}, z^*_{y_i}; y_i)\bigr)}{\| D(\tilde z^k_{y_i} - z^*_{y_i}) \|^2} \;\approx\; 3.11
\end{equation}
is the largest constant for which \cref{assumption:quadraticgrowth} holds across every retained iterate; the dashed line in \cref{fig:qgassumption} traces the corresponding upper bound $2(F - F^*)/\hat\mu$ on $\| D(z - z^*)\|^2$, and every iterate sits on or below it.
Since $\hat\mu$ is a minimum over finitely many iterates and uses one computed minimizer $z^*_y$, it is an empirical estimate and not a certified bound on the constant in \cref{assumption:quadraticgrowth}.
The true constant can be larger leading to faster convergence.

\begin{figure}[htbp]
    \centering
    \includegraphics{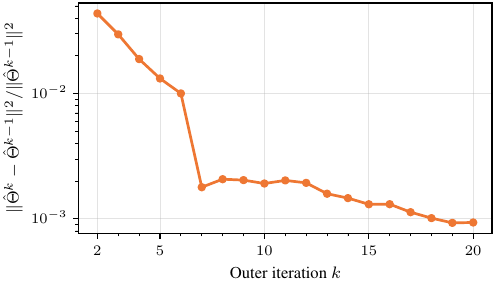}
    \caption{Relative parameter drift of the learned model across outer iterations $k$. We can see a decrease of the relative drift across the training epochs even though the training steps remain constant for later iterations.}
    \label{fig:drift}
\end{figure}

\myparagraph{Bounded parameter drift}
The complexity term in \cref{thm:pacbayes} depends on the cumulative parameter drift $\sum_k \|\hat\Theta^k - \hat\Theta^{k-1}\|^2$ through $\mathrm{KL}(\mathbb Q \| \mathbb P)$, where $\hat\Theta^k$ stacks the trained weights of both subproblem networks at iteration $k$.
\Cref{fig:drift} reports the relative drift $\|\hat\Theta^k - \hat\Theta^{k-1}\|^2 / \|\hat\Theta^{k-1}\|^2$ produced by greedy training with warm-starting over $20$ outer iterations.
The drift decays by roughly an order of magnitude over the first six iterations, drops sharply once the per-iteration training budget is reduced, and continues to shrink towards $10^{-3}$ thereafter.

We observe this decay empirically.
A plausible explanation is that consecutive subproblems become similar as the scheme converges, so warm-started fine-tuning needs only small parameter adjustments.
Since the late increments are small, additional trained iterations add little to the divergence, so the complexity term grows only slowly with the trained depth.

\begin{figure}[htbp]
    \centering
    \includegraphics[width=\linewidth]{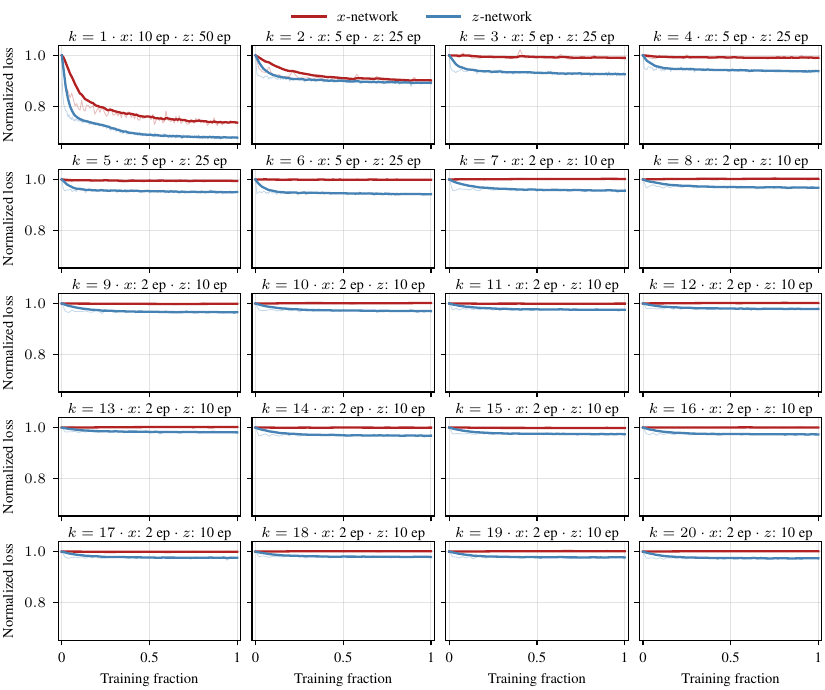}
    \caption{Per-iteration training loss of the $x$- and $z$-networks under greedy training with warm-starting, normalized to the loss at the start of each outer iteration. Panel headers list the epoch budget used for each network at iteration $k$.}
    \label{fig:trainloss}
\end{figure}

\myparagraph{Efficient warm-starting of the training}
Greedy training (\cref{alg:layerwisetraining}) initializes the networks at iteration $k$ from the parameters trained at iteration $k-1$.
\Cref{fig:trainloss} reports the per-outer-iteration training curves of the $x$- and $z$-networks, each normalized to its loss at the start of the iteration, where the panel headers show the epoch budget allocated to each network.
The first iteration is trained from cold-initialized weights and consumes the largest budget, with both networks reducing their loss substantially.
From the third iteration the $x$-network is essentially at a fixed point, while the $z$-network continues to exhibit only mild residual fitting.
The flattening of both curves justifies the budget cuts visible in the panel headers and supports the claim of \cref{sec:greedytr} that warm-starting reduces the per-iteration training cost at later iterations.

\subsection{Learned method comparison}
\begin{figure}[htbp]
    \centering
    \includegraphics[width=\linewidth]{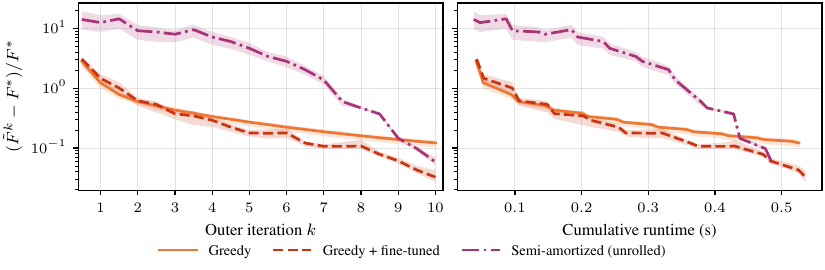}
    \caption{Comparison of three training strategies for the learned alternating scheme at $K=10$ outer iterations: greedy training, greedy training with end-to-end fine-tuning, and the semi-amortized variant. Both panels show the relative objective gap $(\tilde{F}^K - F^*)/F^*$ on a log scale, against the outer iteration $k$ on the left and against cumulative runtime on the right. The shaded bands cover the $5$th--$95$th percentile range over the test instances.}
    \label{fig:learnedmethods}
\end{figure}

In \cref{fig:learnedmethods} we compare three training strategies that share the architecture introduced in \cref{sec:learneddetails} but differ in how the per-iteration networks are fitted: (i) \emph{greedy training} via our proposed \cref{alg:layerwisetraining}; (ii) \emph{greedy training with end-to-end fine-tuning}, where the $K$ greedily trained network pairs are subsequently fine-tuned jointly through the full alternating scheme; and (iii) a \emph{semi-amortized} variant that retains the per-iteration architecture of the proposed scheme but trains all $K$ network pairs jointly through the unrolled alternating scheme.

Greedy training reaches a slightly higher final-iteration gap than the unrolled baseline, consistent with the objective values reported in \cref{sec:learnedcomparison}.
Adding end-to-end fine-tuning to the greedy networks closes this gap and produces the lowest final-iteration error of the three variants.
Interestingly, the fine-tuned variant is able to reduce the functional value further than the greedy variant across most iterations.
In contrast, the unrolled variant starts at a markedly higher relative gap and decreases monotonically as $k$ approaches $K$, since end-to-end training penalizes only the final iterate and does not train on the intermediate iterates of the alternating scheme.

\subsection{Comparison to learned warm-starting}
Using learned warm-starting methods \citep{sambharya2024learning} is an orthogonal idea to our acceleration of the individual subproblems.
The original method is warm-started using the FBP solution as we describe in \cref{subsec:problem}.
In the experiment here, we instead reuse the solution obtained from the fully-amortized model as a warm start to the alternating minimization scheme.

\begin{figure}[htbp]
    \centering
    \includegraphics{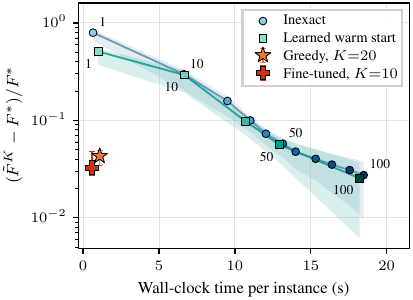}
    \caption{Pareto-front comparison of the learned warm start with the inexact and learned methods. Both iterative methods are run for 20 outer iterations with different number of inner iterations. The inner solver can be stopped early if the subproblem has converged.
    Markers show the median over test instances and bands the 5th--95th percentiles and numbers give the inner-iteration budget.}
    \label{fig:warmstars}
\end{figure}

The results are shown in \cref{fig:warmstars} where we can see that the model performs very similar to the inexact method and only when restricting the methods to a single inner iteration there is a benefit to the learned warm starting.
Both methods are dominated by the learned approach which leads to significantly better performance.
For the single iteration problem the learned method is slightly slower overall due to the additional computation required to predict the warm-start.
When using additional iterations for the inner solvers some of this benefit vanishes since the x-subproblems are faster to converge for the warm-started version leading to a slight overall benefit in the subproblem computations.
However, overall both methods perform nearly on par when using sufficiently many inner iterations.

This is consistent with our theory: in the bound of \cref{theorem:linearbound}, the warm start enters explicitly only through the initial Lyapunov energy, which is discounted exponentially in later iterations, although it also changes the later errors.
Thus, this does not effect the overall convergence speed as much as the later error terms.
Instead, the convergence is dominated by the second term in the bound which is the suboptimality of the inner iterations that can lead to a large gap to the true function value.

\FloatBarrier
\subsection{Evaluation of the PAC-Bayes bound}
\label{app:pacbayes-ct}
\begin{table}[b]
    \centering
    \caption{Evaluation of \cref{eq:pacbayes} for CT reconstruction ($K = 20$, $N = 2\thinspace000$, $\rho = 0.243$) over the prior standard deviation $\varsigma$.
    The risk term is $\sum_k \rho^{K-k}\mathbb E_\mathbb Q[\hat\Delta^k]$, and the trivial bound grows with $\varsigma$ because the plug-in estimate $\hat M$ does.}
    \label{tab:pacbayes-ct}
    \small
    \begin{tabular}{lrrrrrr}
    \toprule
    variance $\varsigma$ & $\mathrm{KL}$ & risk term & complexity & RHS & trivial & RHS/trivial \\
    \midrule
    $0.001$ & $1.9 \times 10^{9}$ & $2.2 \times 10^{-3}$ & $608$ & $608$ & $0.83$ & $733$ \\
    $0.0032$ & $1.9 \times 10^{8}$ & $0.080$ & $192$ & $192$ & $0.83$ & $232$ \\
    $0.01$ & $1.9 \times 10^{7}$ & $1.57$ & $214$ & $216$ & $2.92$ & $74$ \\
    $0.032$ & $1.9 \times 10^{6}$ & $4.11$ & $139$ & $143$ & $5.99$ & $24$ \\
    $0.1$ & $1.9 \times 10^{5}$ & $6.45$ & $109$ & $116$ & $14.9$ & $7.8$ \\
    $0.32$ & $1.9 \times 10^{4}$ & $38.7$ & $237$ & $276$ & $102$ & $2.7$ \\
    $1$ & $1.9 \times 10^{3}$ & $1.1 \times 10^{3}$ & $7.1 \times 10^{3}$ & $8.2 \times 10^{3}$ & $9.7 \times 10^{3}$ & $0.85$ \\
    \bottomrule
    \end{tabular}
\end{table}

We evaluate the bound with the procedure of \cref{app:pacbayes-eval}.
The CT data used for the training does not satisfy the assumptions in \cref{thm:pacbayes}.
The training instances are axial slices from $9$ patients, so neighboring slices of the same patient are correlated and the $N$ samples are not i.i.d.\ draws from a population distribution.
The bound is therefore not a guarantee even with certified constants, which does not change the conclusion below since it is uninformative for every chosen variance $\varsigma$.

We sample $5$ posterior trajectories per value of $\varsigma$ on $10$ training images.
The grid consists of $200$ equally spaced values of $\alpha$ in $[10^{-4}, 0.999\,\alpha_\mathrm{max}]$, where $\alpha_\mathrm{max}$ is the largest $\alpha$ allowed by \cref{eq:alphabounds}, and the $15$ sampled values of $\varsigma$ between $10^{-6}$ and $10$, so that $G = 3\thinspace000$ and $C = 15.4$.

\myparagraph{Results}
We evaluate the greedy model with $K = 20$ on a fixed subset of $N = 2\thinspace000$ training slices (\cref{subsec:problem}), with $\hat\mu = 3.11$.
\Cref{tab:pacbayes-ct} shows that the bound is uninformative for every $\varsigma$: it exceeds the trivial bound for $\varsigma \leq 0.32$, and for $\varsigma \geq 1$ it lies below the trivial bound only because both exceed the observed discounted gap by more than six orders of magnitude.
The Monte Carlo standard error of the risk term over the $5$ posterior trajectories is at most $14\%$ of its value for $\varsigma \leq 0.1$ and $20\%$ at $\varsigma = 0.32$.
At $\varsigma = 1$ it grows to $69\%$, so within two standard errors the ratio of the bound to the trivial bound lies between $0.69$ and $1.0$.
This error does not account for evaluating the posterior on $10$ of the $N$ training slices, and neither source changes the conclusion.

The main cause for the loose bound is the large parameter drift $D = \sum_k \|\hat\Theta^k - \hat\Theta^{k-1}\|^2 = 3\thinspace819$ of the CT networks.
The ratio between the complexity term and the trivial bound is approximately $\tfrac34 \sqrt{(\mathrm{KL} + C)/N}$, which does not depend on $M$.
It falls below one only if $\mathrm{KL} < 16N/9 - C$, that is, only if $\varsigma \gtrsim 0.73$.
Gaussian perturbations of this size destroy the networks: at $\varsigma = 1$, the posterior risk exceeds $10^3$.
Conversely, for $\varsigma \leq 3 \times 10^{-4}$ the posterior risk matches that of the deployed networks, but then the divergence exceeds $10^{10}$.

The risk side of the bound, in contrast, is small and generalizes.
At the deployed weights, the discounted training risk is $1.13 \times 10^{-3}$.
The per-iteration training risk $\hat\Delta^k$ decreases from $6.5 \times 10^{-2}$ at $k = 1$ to $8.2 \times 10^{-4}$ at $k = 20$, and the test risk follows the same profile.
The obstruction is therefore the size of the networks relative to $N$, measured by the drift, and not the generalization of the learned solver itself.
For a smaller learned update, the bound becomes informative (\cref{app:pacbayes-rc}).

\FloatBarrier
\section{Robust matrix completion}

\subsection{Proof of \cref{prop:qg2}}
\label{app:qg2}
Here, we show the proof for the quadratic growth condition for the robust matrix completion problem.
Our approach follows the standard dual certificate analysis of robust PCA \citep{candes2011robust, chandrasekaran2011} and the quadratic growth literature for structured convex problems \citep{cui2017, drusvyatskiy2018error, zhou2017unified}.

\input{proofs/proposition_2}

\subsection{Verification of the quadratic growth condition}
\label{app:qg2-verification}
The constant $\mu$ from the proof certifies existence but is conservative in practice.
We solve all 128 test instances to a KKT residual of at most $1.4 \times 10^{-12}$ and compute the margins at the resulting optimum.
Both regularity conditions hold in practice on every instance, with $\delta_1 \in [8.4 \times 10^{-4}, 0.11]$, $\delta_2 \in [3.3 \times 10^{-5}, 0.026]$, and $\sigma \in [0.050, 0.27]$.
By \cref{prop:qg2}, the minimizer is therefore unique on every test instance.

On every instance, the smallest margin exceeds the KKT residual of that instance by more than three orders of magnitude, so the conditions are resolved numerically.
The small margin $\delta_2$, which comes from residual entries close to the threshold $\lambda_2$, makes the formula give a value of $\mu$ many orders of magnitude below the empirical constant along the iterates.
Along the iterates of proximal gradient and ADMM, the smallest ratio $2(F_y - F^*_y)/\|P_\Omega(z - z^*_y)\|^2$ of each instance lies in $[0.56, 0.77]$.
Along the deployed iterates of the learned models of \cref{sec:matrixlearning,app:pacbayes-rc}, the smallest ratio over all $128$ instances lies in $[0.81, 1.90]$, depending on the model.

\subsection{A common constant over instances}
\label{app:qg2-common}
The constant in \cref{prop:qg2} depends on the instance through the margins $\delta_1$, $\delta_2$, and $\sigma$, while \cref{thm:pacbayes} uses one common constant $\mu$ for all instances.
The two conditions are those under which robust PCA theory establishes exact recovery with high probability for random low-rank and sparse components \citep{candes2011robust, chandrasekaran2011}.
We therefore obtain a common constant on an event of instances and bound the probability of its complement.

Fix margins $\bar\delta_1, \bar\delta_2, \bar\sigma > 0$ before sampling, and let $G$ be the set of data $y$ for which a minimizer of $F_y$ satisfies strict complementarity with $\delta_1 \geq \bar\delta_1$ and $\delta_2 \geq \bar\delta_2$ and restricted injectivity with $\sigma \geq \bar\sigma$.
Since $y$ has bounded support, the iterates are bounded for every parameter value: the learned $x$-update is projected onto a box, and the exact $z$-update satisfies $\|\tilde z\|_\infty \leq \|y\|_\infty + \|\tilde x\|_\infty$ (\cref{sec:matrixlearning}). Hence there is a constant $\bar r$ with $F_y(\tilde x^k_y, \tilde z^k_y) - F^*_y \leq \bar r$ for all $y$, all parameters, and all $k$.

For $y \in G$, the proof of \cref{prop:qg2} with sublevel value $\bar r$ gives quadratic growth along all iterates with a constant $\mu_G > 0$ that depends only on $\bar\delta_1$, $\bar\delta_2$, $\bar\sigma$, $\lambda_1$, $\lambda_2$, and $\bar r$.
Hence \cref{theorem:linearbound} holds for every $y \in G$ with the contraction rate $\rho$ of $\mu_G$.
For $y \notin G$, the resulting gap is at most $\bar r$.
Since the right-hand side of \cref{theorem:linearbound} is nonnegative, we obtain for every $y$
\begin{equation*}
    \tilde F^K_y - F^*_y \leq \rho^K \Phi_0 + \frac{1+\alpha}{\alpha} \sum_{k=1}^K \rho^{K-k} \Delta^k_y + \bar{r} \, \mathds{1} [y \notin G].
\end{equation*}
The proof of \cref{thm:pacbayes} uses \cref{theorem:linearbound} only through this pointwise bound, so \cref{eq:probconv} holds with $\mu = \mu_G$ and the additional $\bar r \, \Pr(y \notin G)$ term.

Membership in $G$ depends only on $y$ and not on the network weights.
If $n$ independent instances all lie in $G$, then $\Pr(y \notin G) \leq \log(1/\beta')/n$ with probability at least $1 - \beta'$, since otherwise all $n$ instances lie in $G$ with probability at most $(1 - \log(1/\beta')/n)^n \leq \beta'$.
A union bound with \cref{thm:pacbayes} gives a total confidence of $1 - \beta - \beta'$.

The margins in \cref{app:qg2-verification} are computed after solving the test instances, so they illustrate the construction rather than certify $G$ for fixed margins.
As discussed, the explicit constant $\mu_G$ is many orders of magnitude below the empirical constant, and we use the empirical constant for the numerical evaluation in \cref{app:pacbayes-rc}.

\subsection{Problem parametrization}
\label{sec:matrixdata}
\myparagraph{Data generation}
Each instance of the dataset is a $100 \times 100$ matrix $y = P_\Omega(L + S)$ with a low-rank component $L$ and a sparse corruption $S$.
The low-rank component is generated as $L = U V^T$ with $U, V \in \mathbb{R}^{100 \times r}$ with rank $r = 5$, and i.i.d.\ entries drawn from $\mathcal{N}(0, 1/100)$.
Each entry of $S$ is nonzero independently with probability $0.05$.
Its nonzero entries are drawn uniformly from the interval $[-3s, 3s]$, where $s$ is the empirical entrywise standard deviation of the low-rank components in the sampled batch, so that the corruption-to-signal ratio does not depend on $n$, $m$, or $r$.
The observation set $\Omega$ contains each entry independently with probability $0.4$, and the observations contain no additional noise.

The low-rank factors are sampled from a Gaussian, so the data $y$ is not bounded by construction.
We therefore define $p_y$ as this generator with the additional constraint that it is conditioned on $\|y\|_\infty \leq c$ with $c = 20\sqrt{r/(nm)} \approx 0.45$, which is fixed from the entry-wise scale of the low-rank component before sampling.
Sampling from the conditioned distribution amounts to rejecting draws with $\|y\|_\infty > c$.
No instance in the training pools or in the test set is rejected in practice and the largest observed value is $\|y\|_\infty = 0.20$, so these samples are exact draws from the conditioned distribution.

\myparagraph{Training and test data}
Since instances are cheap to generate, the networks for \cref{tab:robustcompletion} are trained on a new dataset of instances at every optimization step instead of a fixed training set.
The models used to evaluate the PAC-Bayes bound instead use a fixed training pool of $N$ instances, so that the sample size in \cref{thm:pacbayes} is well defined (\cref{app:pacbayes-rc}).
All methods are evaluated on the same $128$ test instances, generated once with a fixed random seed that is independent of the training seeds.

\myparagraph{Regularization parameters}
We use $\lambda_1 = 0.05$ and $\lambda_2 = 0.01$ for all instances chosen based on cross-validation to minimize the reconstruction error of the minimizer using the ground truth data without corruption.

\myparagraph{Initialization}
All methods in our experiments are initialized with zero vectors $x^0 = z^0 = 0$.

\subsection{Learned method}
\label{sec:matrixlearning}
\myparagraph{Architecture}
The learned scheme keeps the exact closed-form $z$-update and replaces only the $x$-update which is computationally expensive.
At each outer iteration $k$, the network $\phi_{\theta^k}$ performs one proximal gradient step on the $x$-subproblem with a learned step size $\alpha_k > 0$ and a learned threshold vector $\tau_k \in \mathbb{R}^{\min(n, m)}_{>0}$, in the style of LISTA \citep{gregor2010learning}:
\begin{equation*}
    \tilde x^{k+1} = \mathrm{SVT}_{\tau_k}\left( \tilde x^k - \alpha_k P_\Omega(\tilde x^k + \tilde z^k - y) \right),
    \qquad
    \tilde z^{k+1} = P_\Omega\, \mathcal{S}_{\lambda_2}(y - \tilde x^{k+1}).
\end{equation*}
Here, $\mathrm{SVT}_{\tau_k}(u) = U \operatorname{diag}\left(\max(\sigma_i - (\tau_k)_i, 0)\right) V^T$ is a spectral shrinkage map, which is the proximal operator of the nuclear norm when all thresholds are equal, requiring the singular value decomposition (SVD) of the input $u = U \operatorname{diag}(\sigma) V^T$ with singular values sorted in decreasing order.
The proximal operator of the $\ell_1$ norm $\mathcal{S}_{\lambda_2}$ is given by the entrywise soft-thresholding and is efficient to compute.

Unlike the classical SVT, in the learned method the threshold may differ across singular values, so that a layer has more flexibility and can shrink the leading singular values less than the tail.
Entries outside the observed set $\Omega$ do not appear in the objective, and the $z$-update sets them to zero.
Since the $z$-update is the exact minimizer, the error $\zeta^k_y$ in \cref{eq:inexact-error} is zero for every $k$.

Both $\alpha_k$ and $\tau_k$ are parameterized through a softplus function to ensure they are positive.
Each outer iteration therefore has $101$ parameters, $808$ in total for $K = 8$, and costs a single SVD, the same as one iteration of each baseline.
Beyond the trained depth, every iteration $k > 8$ reuses the parameters of the last trained iteration.

To keep the iterates bounded for every parameter, the $x$-update is followed by the entrywise projection onto the box $\{x \mid |x_{ij}| \leq c\}$ with the same a priori constant $c = 20\sqrt{r/(nm)} \approx 0.45$ derived based on the  problem generation from \cref{sec:matrixdata}.
The $z$-update needs no projection, since the exact minimizer already satisfies $\|\tilde z^{k+1}\|_\infty \leq \|y\|_\infty + c$, and projecting it would make the $z$-step inexact.

\myparagraph{Training settings}
We train the networks greedily in a layer-to-layer fashion as in \cref{alg:layerwisetraining}.
The training loss at iteration $k$ is exactly the $x$-subproblem objective $\lambda_1 \|\tilde x^{k+1}\|_* + \frac12 \|P_\Omega(y - \tilde z^k - \tilde x^{k+1})\|_F^2$, averaged over the batch, so training never uses the ground truth $L$ and remains fully unsupervised.

The first iteration is initialized with $\alpha_1 = 1$ and $\tau_1 = \lambda_1 \mathbf{1}$, which reproduces one classical proximal gradient step with step size $1/L_x = 1$.
Every later iteration is initialized from the parameters of the previous network using the greedily trained parameters as a warm start.
Each iteration is trained for $400$ Adam steps with learning rate $10^{-2}$ and batch size $16$, so every iteration sees $6\thinspace400$ fresh instances.
We repeat the training with $5$ random seeds and report the mean over the seeds.

\subsection{Baseline methods}
\label{sec:matrixbaseline}
All baselines run in the same implementation and compute exactly one SVD per outer iteration, so that the outer iteration $k$ is a common budget for all methods in \cref{tab:robustcompletion}.

\myparagraph{FISTA}
The classical alternating scheme solves the $x$-subproblem with FISTA \citep{beck2009fast} with step size $1/L_x = 1$, since $\|P_\Omega\| = 1$.
We apply warm-starting at the previous iterate, and computes the $z$-update in closed form.
In \cref{tab:robustcompletion} we truncate the FISTA iteration to use one inner iteration per outer iteration in  order to match the budget of the other methods.
With a single inner iteration there is no momentum step, so this column coincides with the proximal gradient method discussed below at the fixed step size $\alpha = 1$.
For the time-to-match comparison in \cref{tab:robustcompletiontiming}, we increase the number of inner iterations over the grid $\{1, 2, 3, 4, 5, 6, 8, 10, 15, 20, 30, 50\}$ and report the smallest budget at which FISTA reaches the objective of the learned scheme.

\myparagraph{Tuned proximal gradient}
This proximal-gradient baseline applies the same update as the learned network, but with a fixed step size $\alpha$ and the uniform threshold $\alpha \lambda_1$ in every iteration.
Comparing with a tuned step-size separates the gain from learning the update from the gain of a larger step size itself.
We select $\alpha$ from the grid $\{1, 1.5, 1.8, 2, 2.2, 2.4, 2.5, 2.6, 2.8, 3\}$ and obtain the best performance when choosing $\alpha = 2.6$.
For $\alpha > 2/L_x$, the method has no convergence guarantee, so this is an empirical ablation.

\myparagraph{ADMM}
The problem \eqref{eq:objectiverobust} can be reformulated such that each iteration can be efficiently computed using ADMM \citep{boyd2011distributed} which introduces an auxiliary copy $w$ of $x$ and moves the nuclear norm on it while adding an equality constraint such that this reformulation leads to the equivalent problem
\begin{equation*}
    \begin{array}{ll}
        \text{minimize} & \lambda_1 \|w\|_* + \frac12 \|P_\Omega(y - x - z)\|_F^2 + \lambda_2 \|z\|_1\\
        \text{subject to} & x = w,
    \end{array}
\end{equation*}
where $x$, $z$, and $w$ are the optimization variables.
With a chosen penalty parameter $\rho > 0$ and scaled dual variable $u$, the augmented Lagrangian is given by
\begin{equation*}
    \mathcal{L}_\rho(x, z, w, u) = \lambda_1 \|w\|_* + \frac12 \|P_\Omega(y - x - z)\|_F^2 + \lambda_2 \|z\|_1 + \frac{\rho}{2} \|x - w + u\|_F^2 - \frac{\rho}{2} \|u\|_F^2.
\end{equation*}
One ADMM iteration minimizes $\mathcal{L}_\rho$ over $x$, $z$, and $w$ in turn which can all be solved in closed form and then updates the dual variable:
\begin{align*}
    x^{k+1} &= \left( P_\Omega(y - z^k) + \rho (w^k - u^k) \right) \oslash \left( P_\Omega \mathbf{1} + \rho \mathbf{1} \right), \\
    z^{k+1} &= P_\Omega\, \mathcal{S}_{\lambda_2}(y - x^{k+1}), \\
    w^{k+1} &= \mathrm{SVT}_{\lambda_1 / \rho}(x^{k+1} + u^k), \\
    u^{k+1} &= u^k + x^{k+1} - w^{k+1},
\end{align*}
where $\oslash$ denotes entry-wise division and $\mathbf{1}$ is the all-ones matrix.
The $x$-update is the closed-form minimizer of a separable quadratic form, the $z$-update is the same soft-thresholding step as in the baseline alternating scheme, and the $w$-update is the proximal operator of $(\lambda_1/\rho)\|\cdot\|_*$.
Only the update over the auxiliary variable $w$ requires an SVD, so one ADMM iteration costs one SVD, the same as one iteration of the learned scheme and the other baseline methods in the comparison.

All variables are initialized at zero.
Since the variable $z$ does not enter the constraint, the $x$- and $z$-updates form one pass of block minimization rather than a joint minimization over $(x, z)$.
We therefore use ADMM as an empirical baseline and do not rely on its convergence theory.
We select the penalty parameter $\rho$ from $\{0.02, 0.05, 0.1, 0.2, 0.3, 0.5, 1, 2, 5\}$ and obtain the best performance for $\rho = 0.2$.
The two copies $x$ and $w$ agree only in the limit, and intermediate iterates of $x$ are generically full rank.
We therefore report the objective value at whichever of $(x, z)$ and $(w, \mathcal{S}(w))$ attains the lower value, where $\mathcal{S}(w)$ is the exact $z$-update at $w$.

\myparagraph{Hyperparameter selection}
Both $\alpha$ and $\rho$ in the proximal gradient and ADMM method respectively are selected by the lowest objective after $8$ outer iterations which matches the trained depth of the learned scheme.
We then use the fixed found parameters for all steps in the outer scheme also beyond the tuning horizon.
Both selected values lie in the interior of their grids.

\myparagraph{Reference solution}
The relative optimality gap in \cref{tab:robustcompletion} is computed as $(F - F^*)/F^*$, where $F$ and $F^*$ are the objective values averaged over the $128$ test instances.
We compute $F^*$ with the classical alternating scheme with $30$ warm-started FISTA iterations per subproblem and $60$ outer iterations, that is, $1\thinspace800$ SVDs per instance.
The averaged objective agrees to all reported digits after $40$, $50$, and $60$ outer iterations.
Since the minimizer is unique on every test instance (\cref{app:qg2-verification}), the minimum $F^*$ does not depend on the choice of solver.

As in the CT case, the reference solution does not match the ground truth: the relative error $\|x^* - L\|_F / \|L\|_F$ is $0.195$, which is the recovery error of the convex formulation itself.

\subsection{Computational resources}
\label{sec:matrixresources}
All robust matrix completion experiments run on an 11th Gen Intel Core i7-11850H CPU (8 cores, 16 threads, 2.50~GHz) and $32$~GB of RAM, without GPU acceleration.
We use PyTorch with its default CPU backend for the SVD.
Training one model takes about $3$ minutes, so the $5$ seeds take about $15$ minutes in total.
Each PAC-Bayes evaluation in \cref{app:pacbayes-rc}, including training on the fixed pool, the reference solves, and the posterior Monte Carlo sampling, takes about $30$ to $40$ minutes on the same machine, and we run it for $5$ different configurations.

\subsection{Additional results}
\label{sec:matrixadditional}

Here, we present additional results on timing, detailed optimality gaps and analyze the learned parameters of the scheme.

\input{tables/robust_completion_timing}

\myparagraph{Timing results}
For the wall-clock times in \cref{tab:robustcompletiontiming}, each method processes all $128$ test instances as one batch.
We perform one SVD before timing so that the initial library initialization does not fall into the first timed method.
The objective itself is evaluated outside the timed region, since evaluating the nuclear norm costs another SVD.
For the learned scheme, we report the mean and standard deviation over the $5$ training seeds.

The SVD dominates the cost of every method.
One outer iteration takes about $0.11$~s for the $128$ instances, both for the learned scheme and for FISTA with one inner iteration.
The learned update adds only a scalar multiplication with the learned step size and a thresholding step with a vector of $100$ entries, so its overhead is negligible which is confirmed by the empirical results above.
At equal computational budget, the wall-clock time is therefore proportional to the number of SVDs, and the gains in the objective decrease from \cref{tab:robustcompletion} carry over directly to runtime.

The final time-to-match column in \cref{tab:robustcompletiontiming} shows how much the learned acceleration is worth for this problem.
At the trained depth $k = 8$, FISTA needs $40$ SVDs and $4.4$~s to reach the objective that the learned scheme attains with $8$ SVDs in $0.87$~s, a $5\times$ overall reduction.
The reduction is $3\times$ at $k = 1$ and $5\times$ at $k = 4$.
For comparison, ADMM with the tuned $\rho = 0.2$ needs $17$ iterations and $1.9$~s to match the learned scheme at $k = 8$, which is still a $2.2\times$ reduction.
Beyond the trained depth, the reduction relative to FISTA shrinks to $2.9\times$ at $k = 16$ and $2.0\times$ at $k = 32$.
This is consistent with \cref{tab:robustcompletion}: once the learned scheme reuses the last trained parameters, its gap plateaus, and the classical methods catch up as also predicted by theory.

\myparagraph{Relative optimality gap}
\Cref{tab:robustcompletion} lists the relative objective gap values of \cref{fig:rcconvergence} at five different budgets.
The FISTA method is truncated and uses one inner iteration per outer iteration and therefore coincides with the proximal gradient method at $\alpha = 1/L$ to match the computational budget.
Beyond the trained depth, the gap of the learned scheme plateaus around $4 \times 10^{-5}$ from $k \approx 14$ onward.
In this regime the bound is dominated by the error term of the learned optimizer.

\input{tables/robust_completion_comparison}

\begin{figure}[htbp]
    \centering
    \includegraphics[width=\linewidth]{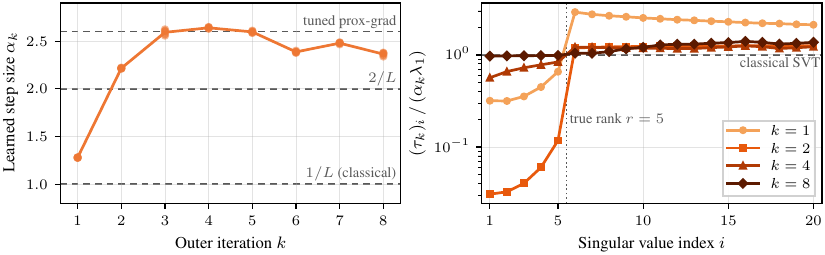}
    \caption{Learned parameters for robust matrix completion solver.
    Left: learned step size $\alpha_k$ for each of the $5$ seeds (markers) and their mean (line).
    Right: learned threshold $(\tau_k)_i$ for the $i$-th singular value relative to the classical threshold $\alpha_k \lambda_1$ at the same step size, averaged over the seeds.
    We show the first $20$ indices only.
    The thresholds of smaller singular values receive no gradient during training, since these singular values lie below the threshold, and, therefore, remain at their initial value.}
    \label{fig:rcparams}
\end{figure}

\myparagraph{Learned parameters}
\Cref{fig:rcparams} shows what the $808$ learned parameters encode.
The step size (left) grows from $1.28$ in the first iteration to about $2.6$ and then settles between $2.35$ and $2.5$, well above the classical $1/L$ and the $2/L$ theoretical limit for the convergence of the scheme.
This roughly matches the chosen step size of the proximal gradient method where the tuned step size $\alpha = 2.6$ is selected for the best performance.
However, the step size alone does not explain the gain in performance, since the proximal gradient method is substantially slower than our learned method as discussed in \cref{tab:robustcompletion}.

The thresholds explain the rest and give the network sufficiently large expressiveness to accelerate the convergence.
In the first iterations, the network shrinks the first $5$ singular values far less than the classical SVT which uses a fixed threshold parameter across all singular values and the remaining ones more.
This is a learned soft rank cutoff at exactly the rank $r = 5$ of the ground truth, learned from the objective alone rather than hand-engineered.
In later iterations, the thresholds approach the classical SVT, as the iterates are already close to the minimizer of the objective.
All learned parameters are nearly identical across the $5$ seeds.

\FloatBarrier
\subsection{Evaluation of the PAC-Bayes bound}
\label{app:pacbayes-rc}
We follow the procedure of \cref{app:pacbayes-eval}.
The bound requires a fixed training sample, so we retrain the learned scheme with $K = 8$ on a fixed pool of $N = 2\thinspace048$ instances with the settings of \cref{sec:matrixlearning}, drawing minibatches from the pool.
These models include the output projection onto the box $\{x \mid |x_{ij}| \leq 0.45\}$ (\cref{sec:matrixdata}).
The projection is inactive on all iterates at the selected $\varsigma$ and becomes active only for $\varsigma \geq 2.3$.
Since the $z$-update is exact, $\Delta^k_y = \xi^k_y$.
We compute the reference minimizer of each $x$-subproblem with $40$ FISTA iterations warm-started at the learned iterate and keep the best objective value, so that $\xi^k_y \geq 0$.
This still remains an approximation of the true error since the reference is not exact but ensures non-negativity of the function gap.

The posterior risks use $4$ sampled trajectories per value of $\varsigma$ on $48$ training instances, and the evaluation of the left-hand side uses $4$ trajectories on $48$ held-out instances.
We sample at $12$ values of $\varsigma$ between $10^{-3}$ and $5$.
We further train variants with an explicit drift penalty $\lambda_\mathrm{D} \|\Theta^k - \Theta^{k-1}\|^2$ added to the greedy loss of iteration $k$, which directly targets the divergence in \cref{thm:pacbayes}.
Since \cref{eq:pacbayesraw} holds for all posteriors simultaneously, comparing these models afterwards requires no further union bound.

We use one constant $\mu = 0.564$ for all models.
It is the smallest quadratic growth ratio on the $128$ test instances over three sets of iterates: the deployed iterates of all ten learned models (\cref{sec:matrixlearning} and the models of this section), the iterates under $\mathbb Q$ for every sampled $\varsigma \leq 0.5$ with $8$ trajectories, and the iterates of proximal gradient and ADMM.
The minimum is attained along the classical iterates.
Along the deployed iterates the ratio is at least $0.81$, and under $\mathbb Q$ it is at least $0.66$.
For $\varsigma \geq 1$, the ratio under $\mathbb Q$ drops to $0.30$, so we restrict the grid to $\varsigma \leq 0.5$.
Since the same $\mu$ applies to all models, the rows of \cref{tab:pacbayes-rc-frontier} share the same contraction rate $\rho$.

We use $60$ equally spaced values of $\alpha$ in $[0.02, 0.98\,\alpha_\mathrm{max}]$ with $\alpha_\mathrm{max} = 0.32$ the largest value allowed by \cref{eq:alphabounds} for $\mu = 0.564$, and $300$ log-spaced values of $\varsigma$ in $[10^{-3}, 0.5]$.
Thus $G = 18\thinspace000$ and $C = 17.2$.\looseness=-1

\myparagraph{Results}
\Cref{tab:pacbayes-rc} lists all quantities for the model with $\lambda_\mathrm{D} = 10^{-4}$.
The right-hand side of \cref{eq:pacbayes} is $0.218$, a factor $1.7$ below the trivial bound of $0.375$.
Unlike for CT, the evaluated bound is therefore informative.
It is still about $13$ times larger than the estimated left-hand side of $0.017$, and the complexity term accounts for $92\%$ of it.
The Monte Carlo standard errors of the posterior training risk and of the left-hand side are $2.0 \times 10^{-3}$ and $2.3 \times 10^{-3}$, small compared with the complexity term.\looseness=-1

\input{tables/pacbayes_robust}

The training and test risks at the trained weights agree to within $1\%$ at every iteration, from $8.05 \times 10^{-2}$ and $8.01 \times 10^{-2}$ at $k = 1$ to $6.70 \times 10^{-5}$ and $6.67 \times 10^{-5}$ at $k = 8$.
At the selected $\varsigma$, the posterior risk is larger than the risk at the trained weights in late iterations, {$2.9 \times 10^{-3}$} instead of $6.7 \times 10^{-5}$ at $k = 8$, but the discounted sum grows only by a factor of {$2.0$}.

The full statement \cref{eq:probconv} is vacuous.
For the model with $\lambda_\mathrm{D} = 10^{-4}$ and $\alpha$ selected to minimize \cref{eq:probconv} (yielding $\alpha = 0.14$, $\rho = 0.79$, $(1+\alpha)/\alpha = 8.2$), its right-hand side is $3.13$, while the trivial bound is $4.50$, the initial energy is $\mathbb E_y[\Phi_0] = 0.93$, and the measured gap is $\mathbb E_\mathbb Q \mathbb E_y[\tilde F^8_y - F^*_y] = 5.1 \times 10^{-3}$.
The bound thus exceeds the initial energy.
Most of its right-hand side is the bound on the error sum, $0.37$ at this $\alpha$, multiplied by the factor $(1+\alpha)/\alpha$.

\Cref{tab:pacbayes-rc-frontier} shows how the drift penalty and the sample size change the bound.
The penalty $\lambda_\mathrm{D} = 10^{-4}$ reduces the drift by a factor $2.7$ and tightens the bound from $0.341$ to $0.218$, while the trained model also attains a lower optimality gap.
Larger penalties tighten the bound further, down to $0.095$, but at the cost of a less accurate model.
With only $N = 128$ training instances, the bound is vacuous at the same drift.
A union bound over all horizons $K \leq 32$ (see \cref{app:reuse}) adds $\log 32 \approx 3.5$ to $\mathrm{KL} + C \approx 968$ and leaves these values unchanged at the reported precision.

\input{tables/pacbayes_robust_frontier}

\clearpage
\input{checklist.tex}

\end{document}

%% file: common.tex
\usepackage{graphicx}

\usepackage{amsmath}
\usepackage{amssymb}
\usepackage{amsthm}
\usepackage{dsfont}

\usepackage{bm}

\usepackage{subcaption}
\newcommand{\argmin}{\operatorname*{argmin}}

\newcommand{\Real}{\mathbb{R}}

\usepackage{booktabs}

\usepackage{algorithm}
\usepackage{algpseudocode}

\usepackage{xcolor}
\definecolor{color_cite}{RGB}{5, 50, 200}
\definecolor{color_link}{RGB}{200, 50, 5}

\usepackage{url}
\usepackage[colorlinks=true]{hyperref}
\hypersetup{citecolor=color_cite, linkcolor=color_link, urlcolor=color_link}
\usepackage[capitalise]{cleveref}

\usepackage{subcaption}

\newcommand{\matr}{A}

\usepackage{xcolor}
\usepackage{thmtools}
\usepackage[framemethod=TikZ]{mdframed}
\usepackage{enumitem}
\mdfsetup{
    skipabove=8pt, 
    skipbelow=0pt,
    hidealllines=false,
    innerleftmargin=10pt,
    innerrightmargin=10pt,
    innertopmargin=8pt,
    innerbottommargin=6pt,
    linewidth=1pt,
    roundcorner=3pt 
}

\declaretheoremstyle[
    mdframed={
        backgroundcolor=white, 
        linecolor=black
    }]{thmstyle}

\declaretheoremstyle[
    mdframed={
        backgroundcolor=gray!5, 
        linecolor=black
    }]{assumptionstyle}

\declaretheoremstyle[
        spaceabove=\parskip,   
        headfont=\itshape,
        notefont=\itshape,
        bodyfont=\normalfont
    ]{remarkstyle}
    
\declaretheorem[style=assumptionstyle]{assumption}
\declaretheorem[style=thmstyle]{theorem}
\declaretheorem[style=thmstyle]{proposition}
\declaretheorem[style=thmstyle]{corollary}
\declaretheorem[style=thmstyle]{lemma}
\declaretheorem[style=remarkstyle]{remark}

\usepackage{tikz}
\usetikzlibrary{spy}

\algtext*{EndFor}
\algtext*{EndIf}
\algtext*{EndWhile}
\algtext*{EndProcedure}
\algtext*{EndFunction}

\newcommand{\myparagraph}[1]{\textbf{#1}\hspace{1em}}

\usepackage[toc,page,header]{appendix}
\usepackage{minitoc}  
\doparttoc 
\faketableofcontents 

\usepackage{wrapfig}
\usepackage{placeins}
\usepackage{graphicx}

%% file: algorithms/training.tex
\begin{algorithm}[t]
    \caption{Greedy training algorithm for learned inexact alternating minimization}
    \label{alg:layerwisetraining}
    \begin{algorithmic}[1]
        \State \textbf{input:} Training data $\mathcal{D} = \{ y_i \}_{i=1}^N \subseteq \mathcal{Y}$ of i.i.d. samples, initial parameters $(\theta^0, \vartheta^0)$.
        \For{$k = 0, 1 \dots K-1$}
        \State Warm-start the learning task $(\theta^{k+1}, \vartheta^{k+1}) \leftarrow (\theta^{k}, \vartheta^{k})$
        \State \textit{// Phase 1: train the $x$-network with $\Theta^{1:k}$ held fixed}
        \State Update $\theta^{k+1}$ by SGD to minimize the empirical risk $\frac{1}{N}\sum_{i=1}^{N} F_{y_i}(\tilde{x}^{k+1}_{y_i}(\theta^{k+1}),\, \tilde{z}^k_{y_i})$
        \State \textit{// Phase 2: train the $z$-network with $\Theta^{1:k}$ and $\theta^{k+1}$ held fixed}
        \State Update $\vartheta^{k+1}$ by SGD to minimize the empirical risk $\frac{1}{N}\sum_{i=1}^{N} F_{y_i}(\tilde{x}^{k+1}_{y_i},\, \tilde{z}^{k+1}_{y_i}(\vartheta^{k+1}))$
        \EndFor
        \State \textbf{output:} Learned optimizer for the subproblems of each iteration $\{ \phi_{\theta^k}, \psi_{\vartheta^k} \}_{k=1}^K$
    \end{algorithmic}
\end{algorithm}

%% file: tables/learned_comparison.tex
\begin{table}[t]
    \centering
    \caption{Comparison of different learning based approaches to solve \cref{eq:objective-ctrecon} showing the mean objective and standard deviation across $210$ test instances.}
    \label{tab:learnedcomp}
    \begin{tabular}{lcccc}
    \toprule
    Method & Training time $\downarrow$ & Inference steps $\downarrow$ & Inference time $\downarrow$ & Test objective $\downarrow$ \\
    \midrule
    Fully-amortized & 17h28m & \textbf{1} & $\bf{0.35}$s & $45.547 \pm 1.996$ \\
    Semi-amortized & 34h11m & \underline{10} & \underline{$0.51$s} & $24.425 \pm 0.508$\\
    \textbf{Greedy (ours)} & \textbf{11h21m} & \underline{10} & $0.53$s & $25.824 \pm 0.362$ \\
    \quad + fine tuning & \;\;+5h55m & \underline{10} & $0.54$s & $\bf{23.840} \pm 0.468$ \\
    \textbf{Greedy (ours)} & \underline{16h27m} & 20 & $1.05$s & \underline{$24.043 \pm 0.328$} \\
    \bottomrule
    \end{tabular}
\end{table}

%% file: proofs/lemma_1.tex
\begin{proof}
The proof proceeds in three steps:
(step 1) we derive an inequality for each subproblem, containing an individual error term;
(step 2) we sum the two inequalities, exploit a non-positive residual, and bound the resulting inner product that remains from the error terms;
(step 3) we replace the exact minimizer $z^{k+1}$ on the left-hand side with the inexact iterate $\tilde{z}^{k+1}$.

\textbf{Step 1: Base inequalities.}
Let $z^{k+1} := \arg\min_{z \in \mathcal{Z}} F(\tilde{x}^{k+1}, z)$ be the exact minimizer of the $z$-subproblem given $\tilde{x}^{k+1}$.
By generalized strong convexity with respect to $z$ (see \cref{remark1}) and applying \cref{lemma2} we obtain for any $z \in \mathcal{Z}$:
\begin{equation}\label{eq:alt_ineq22}
    F(\Tilde{x}^{k+1}, z^{k+1} ) + \frac{1}{2}\|B(z - z^{k+1})\|^2 \leq F(\Tilde{x}^{k+1}, z )
\end{equation}
Adding $\zeta^{k+1}$ to both sides and using its definition to simplify the left-hand side gives
\begin{equation}\label{eq:alt_ineq22p}
    \underbrace{F(\Tilde{x}^{k+1}, \tilde{z}^{k+1} )}_{=: \tilde{F}^{k+1}} + \frac{1}{2}\|B(z - z^{k+1})\|^2 \leq F(\Tilde{x}^{k+1}, z) + \zeta^{k+1}.
\end{equation}
For the $x$-update using the same argument as in the previous case, in the exact case we have
\begin{equation}
    F(x^{k+1}, \tilde{z}^k) + \frac{1}{2}\|A(x - x^{k+1})\|^2 \leq F(x, \tilde{z}^k)
\end{equation}
from the generalized strong convexity.
Replacing the exact minimizer $x^{k+1}$ with the approximate iterate $\tilde{x}^{k+1}$ introduces an error $\mathcal{E}^{k+1}(x)$ so that for any $x \in \mathcal{X}$
\begin{equation} \label{eq:alt_ineq1_error}
\begin{split}
    F(\Tilde{x}^{k+1}, \tilde{z}^k ) + \frac{1}{2}\|A(x - \Tilde{x}^{k+1})\|^2
    & \leq F(x, \tilde{z}^{k}) + \mathcal{E}^{k+1}(x),
\end{split}
\end{equation}
where the error is defined as
\begin{equation} \label{eq:error}
    \mathcal{E}^{k+1}(x) := \underbrace{F(\Tilde{x}^{k+1}, \tilde{z}^k ) - F(x^{k+1}, \tilde{z}^k )}_{= \xi^{k+1}}
     + \frac{1}{2}\|A(x - \Tilde{x}^{k+1})\|^2  - \frac{1}{2}\|A(x - x^{k+1})\|^2.
\end{equation}

\textbf{Step 2: Sum, simplify, and bound the inner product.}
Summing \cref{eq:alt_ineq22p} and \cref{eq:alt_ineq1_error} and using the explicit form of the objective \cref{eq:Objective} gives a non-positive residual term and the error term on the right-hand side:
\begin{equation} \label{eq:alt_ineq_sum}
    \tilde{F}^{k+1} + \frac{1}{2}\|B(z - z^{k+1})\|^2 \leq F(x,z) + \frac{1}{2}\|B(z - \tilde{z}^{k})\|^2 + \tilde{\mathcal{R}}^{k+1}(x, z) + \mathcal{E}^{k+1}(x) + \zeta^{k+1},
\end{equation}
where the residual term simplifies using the explicit form of the objective \cref{eq:Objective}:
\begin{equation}
\label{eq:inresidualnonneg}
\begin{split}
    \tilde{\mathcal{R}}^{k+1}(x, z) &= \frac{1}{2} \|Ax + B\tilde{z}^k\|^2  +\frac{1}{2} \|A\tilde{x}^{k+1} + B z\|^2 - \frac{1}{2}\|A\tilde{x}^{k+1} +  B\tilde{z}^k \|^2 \\
    &\qquad - \frac{1}{2}\|A x + Bz\|^2- \frac{1}{2}\|A(x - \tilde{x}^{k+1})\|^2 - \frac{1}{2}\|B(z - \tilde{z}^{k})\|^2 \\
    &= -\frac{1}{2} \|A(\tilde{x}^{k+1} - x) + B (\tilde{z}^{k} - z )\|^2 \\
    & = -\frac{1}{2}\|A(\tilde{x}^{k+1} - x)\|^2 + \langle A(x - \tilde{x}^{k+1}),\, B(\tilde{z}^k - z)\rangle - \frac{1}{2}\|B(z - \tilde{z}^k)\|^2.
\end{split} 
\end{equation}
We now use $\tilde{\mathcal{R}}^{k+1}(x, z)$ to simplify the right-hand side of \cref{eq:alt_ineq_sum}.
Substituting the definition \cref{eq:error} of $\mathcal{E}^{k+1}(x)$, two cancellations occur:
\begin{itemize}
    \item first, the terms $\pm\frac{1}{2}\|A(x - \tilde{x}^{k+1})\|^2$ from $\tilde{\mathcal{R}}^{k+1}$ and $\mathcal{E}^{k+1}$ cancel
    \item second, the $-\frac{1}{2}\|B(z - \tilde{z}^k)\|^2$ from $\tilde{\mathcal{R}}^{k+1}$ cancels with the explicit $+\frac{1}{2}\|B(z - \tilde{z}^k)\|^2$ on the right-hand side of \cref{eq:alt_ineq_sum}.
\end{itemize}
This gives a simplified bound where the errors of the two subproblems can be added up to each other:
\begin{equation} \label{eq:after_cancel}
    \tilde{F}^{k+1} + \frac{1}{2}\|B(z - z^{k+1})\|^2
    \leq F(x,z) + \langle A(x - \tilde{x}^{k+1}),\, B(\tilde{z}^k - z)\rangle - \frac{1}{2}\|A(x - x^{k+1})\|^2 + \Delta^{k+1}
\end{equation}
We bound the inner product by decomposing $x - \tilde{x}^{k+1} = (x - x^{k+1}) + (x^{k+1} - \tilde{x}^{k+1})$.
Using the linearity of the inner product we proceed to bound each part individually as follows.

For the first part, the identity
\begin{equation}
\langle a, b\rangle - \frac{1}{2}\|a\|^2 = \frac{1}{2}\|b\|^2 - \frac{1}{2}\|a - b\|^2
\end{equation}
with $a = A(x - x^{k+1})$ and $b = B(\tilde{z}^k - z)$ gives
\begin{equation} \label{eq:part1_bound}
    \langle A(x - x^{k+1}),\, B(\tilde{z}^k - z)\rangle - \frac{1}{2}\|A(x - x^{k+1})\|^2
    \leq \frac{1}{2}\|B(z - \tilde{z}^k)\|^2.
\end{equation}

For the second part, Young's inequality with parameter $\alpha > 0$ and the bound
\begin{equation}
\frac{1}{2}\|A(x^{k+1} - \tilde{x}^{k+1})\|^2 \leq \xi^{k+1}
\end{equation}
which follows from \cref{eq:gen_strong_conv_min} in \cref{lemma2} applied to the minimizer $x^{k+1}$ of $x \mapsto F(x, \tilde{z}^k)$, using \cref{remark1} give
\begin{equation} \label{eq:part2_bound}
\begin{split}
    \langle A(x^{k+1} - \tilde{x}^{k+1}),\, B(\tilde{z}^k - z)\rangle
    &\leq \frac{\alpha}{2}\|B(z - \tilde{z}^k)\|^2 + \frac{1}{2\alpha}\|A(x^{k+1} - \tilde{x}^{k+1})\|^2 \\
    &\leq \frac{\alpha}{2}\|B(z - \tilde{z}^k)\|^2 + \frac{1}{\alpha}\xi^{k+1}.
\end{split}
\end{equation}

Combining the steps \cref{eq:after_cancel,eq:part1_bound,eq:part2_bound} yields
\begin{equation} \label{eq:intermediate_lemma}
    F(\tilde{x}^{k+1}, \tilde{z}^{k+1}) + \frac{1}{2}\|B(z - z^{k+1})\|^2
    \leq F(x,z) + \frac{1+\alpha}{2}\|B(z - \tilde{z}^k)\|^2 + \frac{1+\alpha}{\alpha}\xi^{k+1} + \zeta^{k+1}.
\end{equation}

\textbf{Step 3: Bridge exact to inexact $z$-iterate.}
It remains to replace the exact $z^{k+1}$ on the left-hand side with the inexact iterate $\tilde{z}^{k+1}$.
We decompose $B(z - z^{k+1}) = B(z - \tilde{z}^{k+1}) + B(\tilde{z}^{k+1} - z^{k+1})$.

Using Young's inequality on the quadratic expansion gives 
\begin{equation}
\|u + v\|^2 \geq \frac{1}{1+\varepsilon}\|u\|^2 - \frac{1}{\varepsilon}\|v\|^2
\end{equation}
for any $\varepsilon > 0$.
Applying this result with $u = B(z - \tilde{z}^{k+1})$ and $v = B(\tilde{z}^{k+1} - z^{k+1})$ yields\looseness=-1
\begin{equation} \label{eq:bridging}
    \frac{1}{2}\|B(z - z^{k+1})\|^2
    \geq \frac{1}{2(1+\varepsilon)}\|B(z - \tilde{z}^{k+1})\|^2 - \frac{1}{2\varepsilon}\|B(\tilde{z}^{k+1} - z^{k+1})\|^2.
\end{equation}

Since $z^{k+1}$ minimizes $z \mapsto F(\tilde{x}^{k+1}, z)$, generalized strong convexity (\cref{lemma2,remark1}) gives
\begin{equation}
\frac{1}{2}\|B(\tilde{z}^{k+1} - z^{k+1})\|^2 \leq \zeta^{k+1}.
\end{equation}
Substituting into \cref{eq:intermediate_lemma} and combining the $\zeta^{k+1}$ terms on the right-hand side gives:
\begin{equation}
    \tilde{F}^{k+1} + \frac{1}{2(1+\varepsilon)}\|B(z - \tilde{z}^{k+1})\|^2
    \leq F(x,z) + \frac{1+\alpha}{2}\|B(z - \tilde{z}^k)\|^2 + \frac{1+\alpha}{\alpha}\xi^{k+1} + \Bigl(1 + \frac{1}{\varepsilon}\Bigr)\zeta^{k+1}.
\end{equation}
In order to combine the bounds, we are setting $\varepsilon = \alpha$ which gives $1/(1+\varepsilon) = 1/(1+\alpha)$ and $1 + 1/\varepsilon = (1+\alpha)/\alpha$, so the error terms combine into
\begin{equation}
\frac{1+\alpha}{\alpha}(\xi^{k+1} + \zeta^{k+1}) = \frac{1+\alpha}{\alpha}\Delta^{k+1}.
\end{equation}
This concludes the proof.
\end{proof}

%% file: proofs/theorem_1.tex
\begin{proof}
Since \cref{lemma0p} holds for any $(x, z) \in \mathcal{X} \times \mathcal{Z}$, we evaluate \cref{eq:lyapanov} at an arbitrary $(x^*, z^*) \in Z^*$, for which $F(x^*,z^*) = F^*$:
\begin{equation}
        \label{eq:lyapanov-theorem2-pre}
        F(\Tilde{x}^{k+1}, \tilde{z}^{k+1} ) + \frac{1}{2(1+\alpha)}\|B(z^* - \tilde{z}^{k+1})\|^2 \leq F^* + \frac{1+\alpha}{2}\|B(z^* - \tilde{z}^{k})\|^2 + \frac{1+\alpha}{\alpha}\, \Delta^{k+1}.
\end{equation}
This inequality holds for every $(x^*, z^*) \in Z^*$.
On the left-hand side we lower-bound by the infimum,
\begin{equation*}                                                                                                                                                       
      \|B(z^* - \tilde{z}^{k+1})\|^2 \geq \inf_{(x', z') \in Z^*}\|B(z' - \tilde{z}^{k+1})\|^2,                                                                           
\end{equation*}
which is independent of the chosen $(x^*, z^*)$.
Taking the infimum over $(x^*, z^*) \in Z^*$ on the right-hand side then yields
\begin{equation}                                                                                                                                                        
      \label{eq:lyapanov-theorem2}                                                             
      \tilde{F}^{k+1} + \frac{1}{2(1+\alpha)} \inf_{(x', z') \in Z^*}\|B(z' - \tilde{z}^{k+1})\|^2                                                                        
      \leq F^* + \frac{1+\alpha}{2} \inf_{(x', z') \in Z^*} \|B(z' - \tilde{z}^k)\|^2 + \frac{1+\alpha}{\alpha}\,\Delta^{k+1}.                                          
\end{equation}
Defining the Lyapunov energy
\begin{equation}                                                                                                                                                        
      \Phi_k := \tilde{F}^k - F^* + \frac{1}{2(1+\alpha)} \inf_{(x', z') \in Z^*} \|B(z' - \tilde{z}^k)\|^2,                                                           
\end{equation}
we now bound the $B$-norm term on the right-hand side of the inequality using \cref{assumption:quadraticgrowth}.
Using \cref{assumption:quadraticgrowth} gives the bound on the function value gap
\begin{equation}
    \tilde{F}^k - F^* \geq \frac{\mu}{2} \inf_{(x', z') \in Z^*} \|B(z' - \tilde{z}^k)\|^2,
\end{equation}
so that by applying the inequality and combining the two $B$-norm terms we obtain
\begin{equation}
    \Phi_k \geq \frac{\mu(1+\alpha) + 1}{2(1+\alpha)} \inf_{(x', z') \in Z^*} \|B(z' - \tilde{z}^k)\|^2,
\end{equation}
and in particular $\inf_{(x', z') \in Z^*} \|B(z' - \tilde{z}^k)\|^2 \leq \frac{2(1+\alpha)}{\mu(1+\alpha) + 1}\, \Phi_k$.
Substituting this into \cref{eq:lyapanov-theorem2} where we bring $F^*$ to the left-hand side gives the one-step contraction:
\begin{equation}
    \Phi_{k+1} \leq \frac{(1+\alpha)^2}{1+\mu(1+\alpha)}\, \Phi_k + \frac{1+\alpha}{\alpha}\, \Delta^{k+1}.
\end{equation}
Now, we recursively unroll these iterations to obtain for $\rho = {(1+\alpha)^2} / {(1+\mu(1+\alpha))} < 1$:
\begin{equation}
    \Phi_K \leq \rho^K \Phi_0 + \frac{1+\alpha}{\alpha} \sum_{k=1}^K \rho^{K-k} \Delta^k,
\end{equation}
where $\Phi_0 = \tilde{F}^0 - F^* + \frac{1}{2(1+\alpha)} \inf_{(x', z') \in Z^*} \|B(z' - \tilde{z}^{0})\|^2$.

It remains to show that $\rho < 1$ for the choice of $\alpha$ in \cref{eq:alphabounds}.
This is equivalent to the condition
\begin{equation}
    (1 + \alpha)^2 < 1 + \mu(1 + \alpha)
\end{equation}
for the specific contraction rate $\rho$. Now setting $t = 1 + \alpha$, we can see this is a quadratic in $t$.
We find the roots of this quadratic form which yields that the condition is satisfied for
\begin{equation}
    \alpha < \frac{\mu -2 + \sqrt{\mu^2 + 4}}{2}
\end{equation}
since the quadratic $t^2 - \mu t - 1$ is negative between its two roots and the lower root is negative, this is satisfied for any $\alpha$ in the interval provided in \cref{eq:alphabounds}.
This concludes the proof of the bound.
\end{proof}

%% file: proofs/corollary_1.tex
\begin{proof}
    In all three cases, the first term of the bound in \cref{theorem:linearbound} given by
    \begin{equation}
        \rho^K \left( F(\tilde{x}^0, \tilde{z}^0) - F^* + \frac{1}{2(1+\alpha)} \inf_{(x', z') \in Z^*} \|B(z' - \tilde{z}^0)\|^2 \right)
    \end{equation}
    decays geometrically at rate $\rho$, so it suffices
    to bound the error sum
    $S_K = \frac{1+\alpha}{\alpha}\sum_{k=1}^K \rho^{K-k} \Delta^k$ for the different cases.

    \textbf{Case i.} In the case that $q < \rho$, we substitute $\Delta^k \leq Cq^k$ and factor out $\rho^K$ to obtain:
    \begin{equation}
        S_K \leq \frac{C(1+\alpha)}{\alpha}\, \rho^K
        \sum_{k=1}^K \left(\frac{q}{\rho}\right)^k.
    \end{equation}
    The ratio $q/\rho < 1$ and the infinite sum for $K \to \infty$ converges as a geometric sum with limit $(q/\rho) / (1-(q/\rho))$ giving $S_K = \mathcal{O}(\rho^K)$.
    
    If instead $q > \rho$, we mirror the argument from above and factor out $q^K$:
    \begin{equation}
        S_K \leq \frac{C(1+\alpha)}{\alpha} q^{K} \sum_{k=1}^K \rho^{K-k} q^{k-K} \leq \frac{C(1+\alpha)}{\alpha} q^{K} \sum_{j=0}^{K-1} \left( \frac{\rho}{q} \right)^j
    \end{equation}
    where $\rho/q < 1$ and we set $j = K - k$ to obtain a converging sum with the same argument as above yielding a $\mathcal{O}(q^K)$ convergence.

    If $q = \rho$, the sum simplifies to $S_K = \frac{C(1+\alpha)}{\alpha}K \rho^K$.
    The linear factor $K$ reflects a resonance between the contraction and the error decay: the contraction damps past errors at exactly the rate at which new ones are injected, so no part of the history dominates and all $K$ contributions accumulate equally.
    However, the resulting rate $\mathcal{O}(K \rho^K)$ is still dominated by $(\rho')^K$ for any $\rho' \in (\rho, 1)$, so the conclusion below holds in this asymptotic sense at the boundary.
    
    Thus, in all cases $S_K = \mathcal{O}(\max(\rho, q)^K)$ yielding in both cases a linear convergence of the function gap $\tilde{F}^K - F^* \leq \mathcal{O}(\max(\rho, q)^K)$.

    \textbf{Case ii.a.} (Convergence without a rate)\\
    Suppose only that $\Delta^k \to 0$. We now show that $\sum_{k=1}^{K} \rho^{K-k} \Delta^k \to 0$ as $K \to \infty$.
    
    Choosing a fixed $\epsilon > 0$. Since $\Delta^k \to 0$, there exists $\bar k \in \mathbb{N}$ such that $\Delta^k < \epsilon$ for all $k > \bar k$. Splitting the sum at $\bar k$:
    \begin{equation}
        \sum_{k=1}^{K} \rho^{K-k} \Delta^k
        = \sum_{k=1}^{\bar k} \rho^{K-k} \Delta^k
        + \sum_{k=\bar k+1}^{K} \rho^{K-k} \Delta^k.
    \end{equation}
    
    For the first sum, $\rho^{K-k} \leq \rho^{K-\bar k}$ for all $k \leq \bar k$, so
    \begin{equation}
        \sum_{k=1}^{\bar k} \rho^{K-k} \Delta^k
        \leq \rho^{K-\bar k} \sum_{k=1}^{\bar k} \Delta^k,
    \end{equation}
    where $\sum_{k=1}^{\bar k} \Delta^k$ is a fixed finite constant and $\rho^{K-\bar k} \to 0$ as $K \to \infty$.
    
    For the second sum, bounding $\Delta^k < \epsilon$ and extending the geometric series gives
    \begin{equation}
        \sum_{k=\bar k+1}^{K} \rho^{K-k} \Delta^k
        < \epsilon \sum_{k=\bar k+1}^{K} \rho^{K-k}
        \leq \frac{\epsilon}{1-\rho}.
    \end{equation}
    
    Hence $\limsup_{K \to \infty} \sum_{k=1}^{K} \rho^{K-k} \Delta^k \leq \epsilon / (1-\rho)$, and since $\epsilon$ was arbitrary, the sum vanishes.
    
    \textbf{Case ii.b.} (Polynomial decay)\\
    We substitute the polynomial decay assumption into the bound and split the resulting error term into two separate sums similar to the previous step at the threshold $m = \lfloor K/2 \rfloor $.
    
    For $k \leq m$, we have $K - k \geq K/2$, so we can bound the first sum as
    \begin{equation}
        \Sigma_1 = \sum_{k=1}^{m} \rho^{K-k} \frac{C}{k^q}
        \leq C \rho^{K/2} \sum_{k=1}^{m} \frac{1}{k^q}.
    \end{equation}
    The harmonic sum is bounded by $\mathcal{O}(K^{\max(1-q,\, 0)})$ for $q \neq 1$ and $\mathcal{O}(\log K)$ for $q = 1$. In either case, the exponential factor $\rho^{K/2}$ dominates the convergence as $K \to \infty$.
     
    For $k > m$, we have $k \geq K/2$ and thus we have $1/k^q \leq 2^q / K^q$, giving the following inequality by bounding the geometric sum:
    \begin{equation}
        \Sigma_2 = \sum_{k=m+1}^{K} \rho^{K-k} \frac{C}{k^q}
        \leq C \frac{2^q}{K^q} \sum_{j=0}^{K-m-1} \rho^{j}
        \leq C \frac{2^q}{(1-\rho) K^q}.
    \end{equation}
     
    Combining the two previous terms and plugging them back into the expression for $S_K$ yields the following polynomial decay since the convergence is determined by the second part of the sum:
    \begin{equation}
        S_K = \frac{1+\alpha}{\alpha} \sum_{k=1}^{K} \rho^{K-k} \Delta^k
        \leq \frac{C(1+\alpha)}{\alpha} \left( \frac{2^q}{(1-\rho) \, K^q} + \rho^{K/2} \sum_{k=1}^m \frac{1}{k^q} \right)
        = \mathcal{O}(1/K^q),
    \end{equation}
    Adding the geometrically decaying term to the exponentially decaying $\Phi_0$ from the result in \cref{theorem:linearbound} preserves the rate giving a sublinear convergence with a rate determined by the exponent $q$.

    \textbf{Case iii.}
    To start, we split the error sum at the threshold $\bar k$ from the assumption:
    \begin{equation}
        \sum_{k=1}^{K} \rho^{K-k} \Delta^k = \sum_{k=1}^{\bar k} \rho^{K-k} \Delta^k + \sum_{k=\bar k+1}^{K} \rho^{K-k} \Delta^k.
    \end{equation}
    For the first sum, $\rho^{K-k} \leq \rho^{K-\bar k}$ for all $k \leq \bar k$, so
    \begin{equation}
        \sum_{k=1}^{\bar k} \rho^{K-k} \Delta^k \leq \rho^{K-\bar k} \sum_{k=1}^{\bar k} \Delta^k,
    \end{equation}
    where $\sum_{k=1}^{\bar k} \Delta^k$ is a finite constant and $\rho^{K-\bar k} \to 0$ as $K \to \infty$.
    For the second sum, the tail bound $\Delta^k \leq \Delta^\text{tail}$ and the geometric series give
    \begin{equation}
        \sum_{k=\bar k+1}^{K} \rho^{K-k} \Delta^k \leq \Delta^\text{tail} \sum_{k=\bar k+1}^{K} \rho^{K-k} \leq \frac{\Delta^\text{tail}}{1-\rho}.
    \end{equation}
    Combining the two bounds yields $\sum_{k=1}^{K} \rho^{K-k} \Delta^k \leq \rho^{K-\bar k} \sum_{k=1}^{\bar k} \Delta^k + \Delta^\text{tail}/(1-\rho)$ for all $K > \bar k$, and the first term vanishes as $K \to \infty$.
    Taking $\limsup$ on both sides gives
    \begin{equation}
        \limsup_{K \to \infty} S_K \leq \frac{1+\alpha}{\alpha} \cdot \frac{\Delta^\text{tail}}{1-\rho}.
    \end{equation}
    
This finalizes the proof.
\end{proof}

%% file: proofs/theorem_2.tex
\begin{proof}
    We define a joint prior $\mathbb P$ over the parameters $\Theta^{1:K}$ as the Gaussian random walk
\begin{equation}
    \label{eq:prior}
    \Theta^k \mid \Theta^{k-1} \sim \mathcal N(\Theta^{k-1}, \varsigma^2 I) \quad \text{for } k = 1, \ldots, K,
\end{equation}
where $\Theta^0$ is a fixed initialization chosen independently of the training data and $\varsigma^2 > 0$ is the prior variance.
Any data-independent prior yields a valid PAC-Bayes bound.
Here, we adopt the random walk above because its KL divergence to a suitably chosen posterior reduces to the parameter drift across iterations, giving the bound a direct interpretation in terms of warm-starting.

We choose a corresponding posterior $\mathbb Q$ that follows the same random-walk structure of $\mathbb P$ but keeps track of the trained trajectory $\hat\Theta^{1:K}$ rather than staying around the initialization $\Theta^0$.
Under $\mathbb Q$, each step adds the trained parameter change to the previous parameter and the same Gaussian noise as in the prior:
\begin{equation}
    \label{eq:posterior}
    \Theta^k = \Theta^{k-1} + (\hat\Theta^k - \hat\Theta^{k-1}) + \epsilon^k, \quad \epsilon^k \stackrel{\text{i.i.d.}}{\sim} \mathcal N(0, \varsigma^2 I), \quad k = 1, \ldots, K,
\end{equation}
with $\hat\Theta^k$ the parameters returned by greedy training and $\hat\Theta^0 := \Theta^0$ the shared initialization.
The KL divergence between posterior and prior then reduces to the total parameter drift across the trained iterations:
\begin{equation}
    \label{eq:KLidentity}
    \mathrm{KL}(\mathbb Q \| \mathbb P) = \frac{1}{2\varsigma^2}\sum_{k=1}^K \|\hat\Theta^k - \hat\Theta^{k-1}\|^2.
\end{equation}
We view the discounted error sum as a single loss function we aim to bound with the PAC-Bayes argument in the parameters $(\Theta^{1:K}, y)$,
\begin{equation}
    L_y(\Theta^{1:K}) := \sum_{k=1}^K \rho^{K-k}\, \Delta^k_y(\Theta^{1:k}),
\end{equation}
with empirical counterpart on the training samples $\mathcal{D}= \{ y_i \}_{i=1}^N$
\begin{equation}
    \hat L(\Theta^{1:K}) := \frac1N \sum_{i=1}^N L_{y_i}(\Theta^{1:K}) = \frac1N \sum_{i=1}^N \sum_{k=1}^K \rho^{K-k}\, \Delta^k_{y_i}(\Theta^{1:k}) = \sum_{k=1}^K \rho^{K-k}\, \hat\Delta^k(\Theta^{1:k}).
\end{equation}
where the population risk is replaced by its empirical counterpart in \cref{eq:empiricalrisk}.

The bounded-error assumption must hold for every parameter value in the support of $\mathbb{Q}$, not only at $\hat\Theta^{1:K}$, since the bound averages over $\Theta^{1:K} \sim \mathbb Q$.
In our learned schemes, the output projection $\Pi_{[l, u]}$ is part of each network, so the iterates remain in a fixed box for every parameter value, including those drawn from $\mathbb Q$.
Together with the bounded support of $y$, the iterates and the data lie in a compact set on which the continuous objective $F_y$ and its subproblem minima are bounded, which bounds each $\Delta^k_y$ uniformly.
Since $0 \leq \Delta^k_y \leq M$ for every parameter trajectory and $\sum_{k=1}^K \rho^{K-k} \leq 1/(1-\rho)$, the combined loss across the trajectory of iterations for a single instance $L_y$ takes values in $[0, M/(1-\rho)]$.

Applying the PAC-Bayes inequality for bounded losses with a finite grid \citep[Theorem 2.4]{alquier2024user} to $L_{y_1}, \ldots, L_{y_N}$ yields, with probability at least $1-\beta$ over the sample, simultaneously for all $\lambda \in \Lambda$,
  \begin{equation}
  \label{eq:pacbayesraw}
      \mathbb E_\mathbb Q[\mathbb E_y [L_y]]
      \;\leq\; \mathbb E_\mathbb Q[\hat L]
      + \frac{\lambda M^2}{8N(1-\rho)^2}
      + \frac{\mathrm{KL}(\mathbb Q\|\mathbb P) + C}{\lambda},
  \end{equation}
where $C = \log((2J+1)/\beta)$.
The parameter $\lambda$ must be chosen independently of the sample, so we use the grid $\Lambda := \{2^j : j = -J, \ldots, J\}$ with $|\Lambda| = 2J+1$ points, fixed before training.

The right-hand side in \cref{eq:pacbayesraw} has the form $a\lambda + b/\lambda$ with $a := M^2/(8N(1-\rho)^2)$ and $b := \mathrm{KL}(\mathbb Q\|\mathbb P) + C$.
Its unconstrained minimizer $\lambda^\star = \sqrt{b/a} = \tfrac{2(1-\rho)}{M}\sqrt{2N(\mathrm{KL}(\mathbb Q\|\mathbb P) + C)}$ depends on the data through the KL divergence and therefore cannot be used directly.
We call $J$ large enough if $\lambda^\star \in [2^{-J}, 2^J]$, which can be checked after training.
In that case, the grid contains a point $\lambda = t \lambda^\star$ with $t \in [2^{-1/2}, 2^{1/2}]$, at which $a\lambda + b/\lambda = \sqrt{ab}\,(t + 1/t) \leq \tfrac{3}{\sqrt 2}\sqrt{ab}$.
Inserting $a$ and $b$ gives\looseness=-1
\begin{equation}
  \mathbb E_\mathbb Q[\mathbb E_y [L_y]]
  \;\leq\; \mathbb E_\mathbb Q[\hat L]
  + \frac{3M}{4(1-\rho)}\sqrt{\frac{\mathrm{KL}(\mathbb Q\|\mathbb P) + C}{N}}.
\end{equation}
Expanding $L_y$ and $\hat L$ and exchanging the finite sum with $\mathbb E_\mathbb Q$ and $\mathbb E_y$ due to the linearity of the expectation operator yields
\begin{equation}
  \label{eq:pacbayes}
  \sum_{k=1}^K \rho^{K-k}\, \mathbb E_\mathbb Q\bigl[\mathbb E_y[\Delta^k_y]\bigr] \;\leq\;\sum_{k=1}^K \rho^{K-k}\, \mathbb E_\mathbb Q[\hat\Delta^k] + \frac{3M}{4(1-\rho)}\sqrt{\frac{\mathrm{KL}(\mathbb Q \| \mathbb P) + C}{N}}.
\end{equation}
Finally, \cref{eq:expectbound} holds for every fixed parameter trajectory.
All terms are bounded, so taking the expectation over $\Theta^{1:K} \sim \mathbb Q$ and bounding the error sum by \cref{eq:pacbayes} yields \cref{eq:probconv}.

\end{proof}

%% file: proofs/corollary_2.tex
\begin{proof}
    Fix a number of outer iterations at least as large as the number of trained steps $K \geq K_0$.
    Under the weight reuse after the training horizon, the prior and posterior on $\Theta^{1:K}$ are the random walks \cref{eq:prior} and \cref{eq:posterior} on $\Theta^{1:K_0}$, extended by the deterministic map that fixes $\Theta^k = \Theta^{K_0}$ for $k > K_0$.
    This prior remains independent of the data, and the proof of \cref{thm:pacbayes} holds for any data-independent prior.\looseness=-1
    
    Both distributions are images of their marginals on $\Theta^{1:K_0}$ under the same injective map that copies the weights, so $\mathrm{KL}(\mathbb Q\|\mathbb P) = \mathrm{KL}_{K_0}$.
    The bounded-error assumption holds for the reused networks as well, so \cref{eq:probconv} holds at horizon $K$ with probability at least $1-\beta$.

    On this event, we split the discounted training risk at an index $1 \leq m \leq K$.
    The early block satisfies
    \begin{equation*}
        \sum_{k=1}^{m}\rho^{K-k}\,\mathbb E_\mathbb Q[\hat\Delta^k]
        = \rho^{K-m}\sum_{k=1}^{m}\rho^{m-k}\,\mathbb E_\mathbb Q[\hat\Delta^k],
    \end{equation*}
    and the late block satisfies
    \begin{equation*}
        \sum_{k=m+1}^{K}\rho^{K-k}\,\mathbb E_\mathbb Q[\hat\Delta^k]
        \leq \hat\Delta^{> m}\sum_{j=0}^{K-m-1}\rho^{j}
        \leq \frac{\hat\Delta^{> m}}{1-\rho}.
    \end{equation*}
    Writing $\rho^K \mathbb E_y[\Phi_0] = \rho^{K-m}\rho^{m}\,\mathbb E_y[\Phi_0]$ and collecting the terms yields \cref{eq:probconvfinite}.
    The split is applied to the inequality after the event is fixed, so the bound holds for every $m$ on the same event.
\end{proof}

%% file: proofs/proposition_1.tex
\begin{proof}
    Fix a bounded set of data $\mathcal{Y}_0$ and a bounded set $\mathcal{S} \subseteq \mathcal{X} \times \mathcal{Z}$.
    Here, $\mathcal{X}$ and $\mathcal{Z}$ are finite-dimensional and we work with the objective \cref{eq:objective-ctrecon} without the rescaling by $1/(2\lambda_1)$ and convert the constant itself at the end by the same factor.
    The proof has four steps: we eliminate $x$ by partial minimization, write the reduced problem in a form where the data only enters a strongly convex function, apply a quadratic growth result based on Hoffman's error bound, and return to the original variables.

    \textbf{Step 1: reduction to $z$.}
    For a fixed $z$, the $x$-subproblem is strongly convex since its Hessian is $2(T^*T + \lambda_1 I) \succ 0$.
    Its unique minimizer follows from the normal equations as\footnote{In practice it is infeasible to compute the closed form minimizer due to the high dimensionality of the problem.}
    \begin{equation}
        \label{eq:partialmin}
        x^*(z) = (T^*T + \lambda_1I)^{-1}(T^*y+\lambda_1Dz).
    \end{equation}
    Substituting $x^*(z)$ into the objective and simplifying with $(T^*T + \lambda_1 I)^{-1} T^* = T^* (T T^* + \lambda_1 I)^{-1}$ gives the reduced problem
    \begin{equation*}
        \min_z \; \varphi_y(z) := \|\bar{A} z - W y\|^2 + \lambda_2 \|z\|_1,
        \qquad
        W := (I + T T^* / \lambda_1)^{-1/2},
        \quad
        \bar{A} := W T D.
    \end{equation*}
    The operator $\bar{A}$ does not depend on $y$, and $\|W\| \leq 1$.
    Let $Z^*_{\mathrm{red}} := \argmin_z \varphi_y(z)$.
    Every optimum $(x^*, z^*) \in Z^*$ satisfies $x^* = x^*(z^*)$, so $Z^* = \{(x^*(z), z) : z \in Z^*_{\mathrm{red}}\}$, and both problems have the same optimal value $F^*_y$.
    Since $\varphi_y(z) \geq \lambda_2 \|z\|_1$, the function $\varphi_y$ is coercive and $Z^*_{\mathrm{red}}$ is nonempty.

    \textbf{Step 2: polyhedral reformulation.}
    We split $z = u - v$ with $u, v \geq 0$ and consider
    \begin{equation*}
        \min_{u, v \geq 0} \; f_y(u, v) := g_y\bigl(\bar{A}(u - v)\bigr) + \lambda_2 \mathbf{1}^T (u + v),
        \qquad
        g_y(t) := \|t - W y\|^2.
    \end{equation*}
    The function $g_y$ is $2$-strongly convex with a $2$-Lipschitz gradient for every $y$.
    Taking $u$ and $v$ as the positive and negative parts of $z$ shows that the optimal value is $F^*_y$.
    Since $\|u - v\|_1 \leq \mathbf{1}^T (u + v)$, every minimizer $(\bar u, \bar v)$ satisfies $\varphi_y(\bar u - \bar v) \leq f_y(\bar u, \bar v) = F^*_y$, so $\bar u - \bar v \in Z^*_{\mathrm{red}}$.
    The problem has the form $g(A w) + c^T w$ over the polyhedron $\{w : C w \leq d\}$ with $A = \bar{A} [I, -I]$, $c = \lambda_2 \mathbf{1}$, $C = -I$, and $d = 0$.
    For this form, \citet[Section 4.3]{necoara2019linear} show that for every $M > 0$, every feasible $w$ with $f_y(w) - F^*_y \leq M$ satisfies
    \begin{equation*}
        f_y(w) - F^*_y \geq \frac{\kappa_y}{2} \operatorname{dist}(w, W^*_y)^2,
        \qquad
        \kappa_y = \frac{2}{\theta^2 (1 + 2M + 2 c_y^2)},
    \end{equation*}
    where $W^*_y$ is the set of minimizers, $c_y := \|\nabla g_y(\bar{A} z^*)\|$ for any $z^* \in Z^*_{\mathrm{red}}$, and $\theta$ is the Hoffman constant of the linear system with matrices $A$, $c^T$, and $C$.
    The Hoffman constant depends only on these matrices, which do not depend on $y$.

    \textbf{Step 3: bounds that hold for all $y \in \mathcal{Y}_0$.}
    Let $\bar y := \sup_{y \in \mathcal{Y}_0} \|y\|$.
    For $z^* \in Z^*_{\mathrm{red}}$, we have $\|\bar{A} z^* - W y\|^2 \leq \varphi_y(z^*) \leq \varphi_y(0) = \|W y\|^2 \leq \bar y^2$, so $c_y = 2 \|\bar{A} z^* - W y\| \leq 2 \bar y$.
    For $(\tilde x, \tilde z) \in \mathcal{S}$, the value $\varphi_y(\tilde z) - F^*_y \leq \varphi_y(\tilde z)$ is bounded by a constant $\bar M$ that depends only on $\mathcal{S}$, $\bar y$, $\|\bar{A}\|$, and $\lambda_2$.
    With $M = \bar M$, we obtain the constant $\bar\kappa := 2 / (\theta^2 (1 + 2 \bar M + 8 \bar y^2))$, which bounds $\kappa_y$ from below for every $y \in \mathcal{Y}_0$.

    Now take $(\tilde x, \tilde z) \in \mathcal{S}$, let $\tilde w$ consist of the positive and negative parts of $\tilde z$, so that $f_y(\tilde w) = \varphi_y(\tilde z)$, and let $(\bar u, \bar v)$ be the projection of $\tilde w$ onto $W^*_y$.
    Since $\bar u - \bar v \in Z^*_{\mathrm{red}}$ and the map $(u, v) \mapsto u - v$ has norm $\sqrt{2}$,
    \begin{equation*}
        \operatorname{dist}(\tilde z, Z^*_{\mathrm{red}})^2 \leq \|\tilde z - (\bar u - \bar v)\|^2 \leq 2 \operatorname{dist}(\tilde w, W^*_y)^2 \leq \frac{4}{\bar\kappa} \bigl( \varphi_y(\tilde z) - F^*_y \bigr).
    \end{equation*}

    \textbf{Step 4: return to $(x, z)$.}
    Since $x^*(\tilde z)$ minimizes $F(\cdot, \tilde z)$, we have $F(\tilde x, \tilde z) \geq F(x^*(\tilde z), \tilde z) = \varphi_y(\tilde z)$.
    With $B = D$ and the description of $Z^*$ from Step 1,
    \begin{equation*}
        \inf_{(x', z') \in Z^*} \|B(z' - \tilde z)\|^2 = \inf_{z' \in Z^*_{\mathrm{red}}} \|B(z' - \tilde z)\|^2 \leq \|B\|^2 \operatorname{dist}(\tilde z, Z^*_{\mathrm{red}})^2.
    \end{equation*}
    Combining the inequalities yields the quadratic growth inequality at every point of $\mathcal{S}$ with constant $\bar\kappa / (2 \|D\|^2)$ for the unscaled objective and every $y \in \mathcal{Y}_0$.
    For the objective rescaled by $1/(2\lambda_1)$, whose coupling term has the form of \cref{eq:Objective} with $B = D$, this gives \cref{assumption:quadraticgrowth} with $\mu = \bar\kappa / (4 \lambda_1 \|D\|^2)$.
\end{proof}

%% file: tables/hyperparameters.tex
\begin{table}[t]
    \centering
    \caption{Details about the training settings for the different learned methods}
    \label{tab:learneddetails}
    \small
    \setlength{\tabcolsep}{4pt}
    \begin{tabular}{lcccc}
    \toprule
    Hyperparameter & Greedy & Greedy + FT & Semi-amortized & Fully-amortized  \\
    \midrule
    Outer iterations ($K$) & $20$ & $10$ & $10$ & $1$ \\
    Total parameters & $3\,369\,100$ & $1\,684\,550$ & $1\,684\,550$ &  $2\,213\,957$ \\
    Training epochs & \cref{sec:learneddetails} & greedy epochs + $10$ & $50$ & $50$ \\
    Optimizer  & AdamW & AdamW & AdamW & AdamW \\
    Batch size & 4 ($x$) / 8 ($z$) & 2 & 2 & 2 \\
    LR Scheduler & cosine & cosine & cosine & cosine \\
    \bottomrule
    \end{tabular}
\end{table}

%% file: proofs/proposition_2.tex
We first make the two conditions of \cref{prop:qg2} precise.
Let $(x^\ast, z^\ast)$ be a minimizer of the robust completion objective, $R = P_\Omega (y - x^\ast - z^\ast)$ the residual at the optimum, $T$ the tangent space of the rank-$r$ matrix manifold at $x^\ast = U \Sigma V^T$, and $S$ the support of $z^\ast$.
First-order optimality states that $R$ is a subgradient of $\lambda_1 \| \cdot \|_*$ at $x^\ast$ and of $\lambda_2 \| \cdot \|_1$ at $z^\ast$, that is,
\begin{gather*}
    R = \lambda_1 (U V^T + W), \quad P_T W = 0, \quad \|W\|_2 \leq 1,
    \\
    R_{ij} = \lambda_2 \text{sign}(z^\ast_{ij}) \text{ for } (i,j) \in S, \quad |R_{ij}| \leq \lambda_2 \text{ for } (i,j) \notin S,
\end{gather*}
where $\|\cdot\|_2$ denotes the spectral norm.
\begin{itemize}
    \item \emph{Strict complementarity} holds if the two inequalities hold with a margin: there are $\delta_1, \delta_2 \in (0, 1]$ with $\|W\|_2 \leq 1 - \delta_1$ and $|R_{ij}| \leq \lambda_2 (1 - \delta_2)$ for all $(i,j) \notin S$.
    \item \emph{Restricted injectivity} holds if there is $\sigma > 0$ with $\|P_\Omega(u + v)\|_F^2 \geq \sigma^2 (\|u\|_F^2 + \|v\|_F^2)$ for all $u \in T$ and all $v$ supported on $S$.
\end{itemize}
The first condition controls moving away from the minimizer by leaving the low-rank and sparse structure, and the second controls moving along it.
The proof shows that \cref{prop:qg2} holds with an explicit constant $\mu$ that depends only on $\delta_1$, $\delta_2$, $\sigma$, $\lambda_1$, $\lambda_2$, as well as the chosen sublevel value.

\begin{proof}

    We write $\mathcal{A}(x, z) := P_\Omega(x + z)$, so that the objective is $F(x, z) = \frac12 \|P_\Omega y - \mathcal{A}(x, z)\|_F^2 + \lambda_1 \|x\|_* + \lambda_2 \|z\|_1$.
    Fix a sublevel value $\bar\varepsilon > 0$ and a point $(\tilde{x}, \tilde{z})$ with $F(\tilde{x}, \tilde{z}) \leq F^* + \bar\varepsilon$.
    Let $d = (d_x, d_z) := (\tilde{x} - x^\ast, \tilde{z} - z^\ast)$ and $\varepsilon := F(\tilde{x}, \tilde{z}) - F^* \in [0, \bar\varepsilon]$.
    
    The proof has four steps: we split the gap $\varepsilon$ into a curvature term and a convexity gap, then we bound the component of $d$ that leaves the structure $T \times S$ by the convexity gap and the component along $T \times S$ by the curvature term, and finally combine the two bounds.

    \textbf{Step 1: split of the gap.}
    The coupling term is quadratic, so its expansion around $(x^\ast, z^\ast)$ is exact.
    Since $R$ is supported on $\Omega$, we have $\langle R, \mathcal{A} d \rangle = \langle R, d_x \rangle + \langle R, d_z \rangle$, and therefore
    \begin{equation}
        \label{eq:rc-split}
        \varepsilon = \frac12 \|\mathcal{A} d\|_F^2 + D_x + D_z,
    \end{equation}
    where
    \begin{equation*}
        D_x := \lambda_1 \|\tilde{x}\|_* - \lambda_1 \|x^\ast\|_* - \langle R, d_x \rangle,
        \qquad
        D_z := \lambda_2 \|\tilde{z}\|_1 - \lambda_2 \|z^\ast\|_1 - \langle R, d_z \rangle.
    \end{equation*}
    Both $D_x$ and $D_z$ are nonnegative because $R$ is a subgradient of the respective norm.

    \textbf{Step 2: strict complementarity bounds the component leaving the structure.}
    For the $\ell_1$ part, each entry $(i,j) \in S$ contributes a nonnegative term to $D_z$, and each entry $(i,j) \notin S$ contributes $\lambda_2 |(d_z)_{ij}| - R_{ij} (d_z)_{ij} \geq \lambda_2 \delta_2 |(d_z)_{ij}|$.
    Hence $D_z \geq \lambda_2 \delta_2 \|(d_z)_{S^c}\|_1$, where $(d_z)_{S^c}$ is the restriction of $d_z$ to the complement of $S$.

    For the nuclear-norm part, let $U_2 \Sigma_2 V_2^T$ be a compact SVD of $P_T^\perp d_x$ and set $W' := U_2 V_2^T$.
    Since $P_T^\perp d_x \in T^\perp$, the matrix $W'$ lies in $T^\perp$ with $\|W'\|_2 \leq 1$, so $U V^T + W'$ is a subgradient of $\|\cdot\|_*$ at $x^\ast$, and $\langle W', d_x \rangle = \|P_T^\perp d_x\|_*$.
    The subgradient inequality with $U V^T + W'$ and the duality between the spectral and nuclear norms give
    \begin{equation*}
        D_x \geq \lambda_1 \langle W' - W, d_x \rangle
        = \lambda_1 \left( \|P_T^\perp d_x\|_* - \langle W, P_T^\perp d_x \rangle \right)
        \geq \lambda_1 \delta_1 \|P_T^\perp d_x\|_*.
    \end{equation*}

    Let $\eta := (P_T^\perp d_x, (d_z)_{S^c})$ be the component of $d$ orthogonal to $T \times S$ and $\omega := d - \eta \in T \times S$ the component along it.
    With $\delta := \min(\lambda_1 \delta_1, \lambda_2 \delta_2)$, the bounds $\|M\|_F \leq \|M\|_*$ and $\|M\|_F \leq \|M\|_1$, and \cref{eq:rc-split}, we get
    \begin{equation*}
        \|\eta\|_F \leq \|P_T^\perp d_x\|_F + \|(d_z)_{S^c}\|_F \leq \frac{D_x + D_z}{\delta} \leq \frac{\varepsilon}{\delta}.
    \end{equation*}

    \textbf{Step 3: restricted injectivity bounds the component along the structure.}
    By \cref{eq:rc-split}, $\|\mathcal{A} d\|_F \leq \sqrt{2 \varepsilon}$.
    The operator $\mathcal{A}$ has norm at most $\sqrt{2}$, since $\|P_\Omega(u + v)\|_F \leq \|u\|_F + \|v\|_F$.
    Restricted injectivity applied to $\omega \in T \times S$ then gives the bound
    \begin{equation*}
        \sigma \|\omega\|_F \leq \|\mathcal{A} \omega\|_F \leq \|\mathcal{A} d\|_F + \|\mathcal{A} \eta\|_F \leq \sqrt{2 \varepsilon} + \sqrt{2} \frac{\varepsilon}{\delta}.
    \end{equation*}

    \textbf{Step 4: combination of the bounds.}
    The components $\omega$ and $\eta$ are orthogonal.
    Using $(a + b)^2 \leq 2a^2 + 2b^2$ and $\varepsilon \leq \bar\varepsilon$, we obtain
    \begin{equation*}
        \|d\|_F^2 = \|\omega\|_F^2 + \|\eta\|_F^2
        \leq \frac{4 \varepsilon + 4 \varepsilon^2 / \delta^2}{\sigma^2} + \frac{\varepsilon^2}{\delta^2}
        \leq \left[ \frac{4}{\sigma^2} + \left( \frac{4}{\sigma^2} + 1 \right) \frac{\bar\varepsilon}{\delta^2} \right] \varepsilon.
    \end{equation*}
    This is quadratic growth on the sublevel set $\{ F \leq F^* + \bar\varepsilon \}$,
    \begin{equation*}
        \frac{\mu}{2} \|(\tilde{x}, \tilde{z}) - (x^\ast, z^\ast)\|_F^2 \leq F(\tilde{x}, \tilde{z}) - F^*,
        \qquad
        \mu := 2 \left[ \frac{4}{\sigma^2} + \left( \frac{4}{\sigma^2} + 1 \right) \frac{\bar\varepsilon}{\delta^2} \right]^{-1},
    \end{equation*}
    and $\mu$ depends only on $\delta_1$, $\delta_2$, $\sigma$, $\lambda_1$, $\lambda_2$, and $\bar\varepsilon$.

    Applying this inequality to another minimizer gives $\varepsilon = 0$ and hence $d = 0$, so the minimizer is unique.
    Finally, substituting $\bar{z} := z - y$ writes the objective in the form \cref{eq:Objective} with $A = B = P_\Omega$, $f = \lambda_1 \|\cdot\|_*$, and $g(\bar{z}) = \lambda_2 \|\bar{z} + y\|_1$.
    The substitution is a translation and leaves all distances unchanged.
    Since $Z^* = \{(x^\ast, z^\ast)\}$ and $\|P_\Omega\| = 1$, we have $\|P_\Omega(\tilde{z} - z^\ast)\|_F \leq \|d\|_F$, so \cref{assumption:quadraticgrowth} holds with the same $\mu$.
\end{proof}

%% file: tables/robust_completion_timing.tex
\begin{table}[ht]
    \centering
    \caption{Wall-clock time in seconds for the $128$ test instances of robust matrix completion.
    Equal SVDs: FISTA at the same outer iteration $k$.
    Time to match: FISTA with the smallest number of SVDs (in parentheses) that reaches the objective of the learned scheme at $k$.}
    \label{tab:robustcompletiontiming}
    \begin{tabular}{lccc}
    \toprule
    $k$ & \bf{Learned (ours)} & FISTA (equal SVDs) & FISTA (time to match) \\
    \midrule
    1 \textit{(trained)} & $0.112 \pm 0.002$\,s & $0.116$\,s & $0.335$\,s (3 SVDs) \\
    4 \textit{(trained)} & $0.439 \pm 0.003$\,s & $0.442$\,s & $2.141$\,s (20 SVDs) \\
    8 \textit{(trained)} & $0.873 \pm 0.003$\,s & $0.878$\,s & $4.355$\,s (40 SVDs) \\
    \midrule
    16 \textit{(extrapolated)} & $1.746 \pm 0.009$\,s & $1.764$\,s & $5.114$\,s (48 SVDs) \\
    32 \textit{(extrapolated)} & $3.495 \pm 0.012$\,s & $3.460$\,s & $6.935$\,s (64 SVDs) \\
    \bottomrule
    \end{tabular}
\end{table}

%% file: tables/robust_completion_comparison.tex
\begin{table}[htbp]
    \centering
    \caption{Relative optimality gap $(F-F^\star)/F^\star$ for robust matrix completion problem at equal SVD budget $k$ (one SVD per outer iteration for all four methods).
    Learned results are the mean $\pm$ standard deviation over $5$ independent training runs.}
    \label{tab:robustcompletion}
    \resizebox{\textwidth}{!}{%
    \begin{tabular}{lcccc}
    \toprule
    Outer iterations $k$ & \textbf{Learned (ours)} & ADMM ($\rho=0.2$) & Prox-grad ($\alpha=2.6$) & FISTA ($\alpha=1/L$) \\
    \midrule
    1 \textit{(trained)} & $\mathbf{(3.46 \pm 0.01)\times 10^{-1}}$ & $7.6\times 10^{-1}$ & $4.1\times 10^{0}$ & $\underline{5.0\times 10^{-1}}$ \\
    4 \textit{(trained)} & $\mathbf{(1.67 \pm 0.01)\times 10^{-2}}$ & $\underline{7.5\times 10^{-2}}$ & $1.3\times 10^{-1}$ & $2.5\times 10^{-1}$ \\
    8 \textit{(trained)} & $\mathbf{(6.59 \pm 0.11)\times 10^{-4}}$ & $\underline{9.9\times 10^{-3}}$ & $1.4\times 10^{-2}$ & $1.2\times 10^{-1}$ \\
    \midrule
    16 \textit{(extrapolated)} & $\mathbf{(3.66 \pm 0.53)\times 10^{-5}}$ & $7.3\times 10^{-4}$ & $\underline{1.4\times 10^{-4}}$ & $2.7\times 10^{-2}$ \\
    32 \textit{(extrapolated)} & $(4.38 \pm 0.64)\times 10^{-5}$ & $\mathbf{1.4\times 10^{-5}}$ & $\underline{2.3\times 10^{-5}}$ & $6.6\times 10^{-4}$ \\
    \bottomrule
    \end{tabular}
    }
\end{table}

%% file: tables/pacbayes_robust.tex
\begin{table}[t]
    \centering
    \caption{Quantities entering \cref{eq:pacbayes} for robust matrix completion at the selected $(\alpha, \varsigma)$.
    Model trained on a fixed pool of $N = 2\thinspace048$ instances with drift penalty $\lambda_\mathrm{D} = 10^{-4}$.
    The constants $\mu$ and $\hat M$ are plug-in estimates, not certified constants.}
    \label{tab:pacbayes-rc}
    \small
    \begin{tabular}{llr}
    \toprule
    Symbol & Meaning & Value\\
    \midrule
    $K$ & trained depth & $8$ \\
    $N$ & training sample size & $2\thinspace048$ \\
    $\beta$ & confidence level & $0.05$ \\
    $C$ & $\log((2J+1)/\beta) + \log G$ with $J = 40$, $G = 18\thinspace000$ & $17.2$ \\
    $\mu$ & common quadratic growth constant (\textit{estimated}) & $0.564$ \\
    $\rho$ & contraction rate at the selected $\alpha = 0.02$ & $0.66$ \\
    $\hat M$ & largest observed $\Delta^k_y$ under $\mathbb Q$ (\textit{estimated}) & $0.13$ \\
    $D$ & drift $\sum_k \|\hat\Theta^k - \hat\Theta^{k-1}\|^2$ & $20.2$ \\
    $\varsigma$ & selected prior standard deviation & $0.103$ \\
    $\mathrm{KL}$ & $D / (2\varsigma^2)$ & $951$ \\
    \midrule
    $\sum_k \rho^{K-k}\hat\Delta^k$ & training risk at $\hat\Theta$ & $8.9 \times 10^{-3}$ \\
    $\sum_k \rho^{K-k}\mathbb E_y[\Delta^k_y]$ & test risk at $\hat\Theta$ & $8.9 \times 10^{-3}$ \\
    $\sum_k \rho^{K-k}\mathbb E_\mathbb Q[\hat\Delta^k]$ & training risk under $\mathbb Q$ & $1.8 \times 10^{-2}$ \\
    $\frac{3\hat M}{4(1-\rho)}\sqrt{(\mathrm{KL}+C)/N}$ & complexity term & $0.200$ \\
    \midrule
    & right-hand side of \cref{eq:pacbayes} & $\mathbf{0.218}$ \\
    & trivial bound $\hat M (1-\rho^K)/(1-\rho)$ & $0.375$ \\
    & test risk under $\mathbb Q$ (left-hand side, estimated) & $1.7 \times 10^{-2}$ \\
    \bottomrule
    \end{tabular}
\end{table}

%% file: tables/pacbayes_robust_frontier.tex
\begin{table}[h]
    \centering
    \caption{Effect of the drift penalty $\lambda_\mathrm{D}$ and the sample size $N$ on the bound for robust matrix completion.
    RHS is the right-hand side of \cref{eq:pacbayes} at the selected $(\alpha, \varsigma)$, and LHS is its left-hand side estimated on held-out instances.
    The last column is the relative optimality gap of the trained model at $k = 8$, averaged per instance.
    All rows use the common constant $\mu = 0.564$ and $\rho = 0.66$.}
    \label{tab:pacbayes-rc-frontier}
    \small
    \begin{tabular}{rrrrrrrr}
    \toprule
    $N$ & $\lambda_\mathrm{D}$ & $D$ & $\mathrm{KL}$ & RHS & trivial & LHS & gap at $k = 8$ \\
    \midrule
    $2\thinspace048$ & $0$ & $53.6$ & $2\thinspace524$ & $0.341$ & $0.375$ & $1.5 \times 10^{-2}$ & $6.2 \times 10^{-4}$ \\
    \midrule
    $2\thinspace048$ & $10^{-4}$ & $20.2$ & $951$ & $0.218$ & $0.375$ & $1.7 \times 10^{-2}$ & $4.3 \times 10^{-4}$ \\
    $2\thinspace048$ & $3 \times 10^{-4}$ & $8.5$ & $418$ & $0.135$ & $0.335$ & $1.3 \times 10^{-2}$ & $1.0 \times 10^{-3}$ \\
    $2\thinspace048$ & $10^{-3}$ & $3.5$ & $174$ & $0.095$ & $0.335$ & $1.4 \times 10^{-2}$ & $3.1 \times 10^{-3}$ \\
    \midrule
    $128$ & $10^{-4}$ & $19.4$ & $915$ & $0.954$ & $0.446$ & $1.6 \times 10^{-2}$ & $4.3 \times 10^{-4}$ \\
    \bottomrule
    \end{tabular}
\end{table}

%% file: checklist.tex
\section*{NeurIPS Paper Checklist}

\begin{enumerate}

\item {\bf Claims}
    \item[] Question: Do the main claims made in the abstract and introduction accurately reflect the paper's contributions and scope?
    \item[] Answer: \answerYes{} 
    \item[] Justification: The claims of our paper are clearly listed out in the end of the introduction as a bullet point list with links to the exact sections that contain the corresponding theoretical and experimental results.

\item {\bf Limitations}
    \item[] Question: Does the paper discuss the limitations of the work performed by the authors?
    \item[] Answer: \answerYes{} 
    \item[] Justification: There is a designated paragraph for limitations of our approach with a corresponding heading in \cref{sec:conclusion}. The limited empirical scope is discussed in the future work paragraph in the same section.

\item {\bf Theory assumptions and proofs}
    \item[] Question: For each theoretical result, does the paper provide the full set of assumptions and a complete (and correct) proof?
    \item[] Answer: \answerYes{} 
    \item[] Justification: Yes, the proofs of the claims can be found in \cref{sec:proofs,sec:proofs-learned} and the corresponding sections in the appendix for the CT reconstruction and matrix completion problem families.

    \item {\bf Experimental result reproducibility}
    \item[] Question: Does the paper fully disclose all the information needed to reproduce the main experimental results of the paper to the extent that it affects the main claims and/or conclusions of the paper (regardless of whether the code and data are provided or not)?
    \item[] Answer: \answerYes{} 
    \item[] Justification: We provide the code to reproduce all experiments and figures in the supplement material and provide a link to a publicly accessible github repository in the abstract. Further, \cref{sec:learneddetails} contains a description of all hyperparameters and design decisions for the experimental settings. The problem setup and baseline methods are described in \cref{subsec:problem} and \cref{sec:baseline}. Accordingly, the methods for the matrix completion experiments are described in the sections \cref{sec:matrixdata}, \cref{sec:matrixlearning}, and \cref{sec:matrixbaseline}.

\item {\bf Open access to data and code}
    \item[] Question: Does the paper provide open access to the data and code, with sufficient instructions to faithfully reproduce the main experimental results, as described in supplemental material?
    \item[] Answer: \answerNo{} 
    \item[] Justification: The Mayo Low Dose CT Grand Challenge dataset is publicly accessible after registration but is distributed under CC BY-NC 4.0 / NIH Controlled Data Access, so we cannot redistribute it. We release all code, the pre-trained dictionary, and model checkpoints in the supplement and a public github repository linked in the abstract, together with instructions for reproducing the experiments once the dataset has been obtained.

\item {\bf Experimental setting/details}
    \item[] Question: Does the paper specify all the training and test details (e.g., data splits, hyperparameters, how they were chosen, type of optimizer) necessary to understand the results?
    \item[] Answer: \answerYes{} 
    \item[] Justification: The hyperparameters are described in \cref{sec:learneddetails}. Further details can be found in the code.

\item {\bf Experiment statistical significance}
    \item[] Question: Does the paper report error bars suitably and correctly defined or other appropriate information about the statistical significance of the experiments?
    \item[] Answer: \answerYes{} 
    \item[] Justification: Shaded bands in our figures show the $5-95$ percentile across the test instances. Tables report the mean and standard deviations over the test data slices. Further, we retrain the networks across $5$ independent seeds and report the standard deviation of the objective gap in the result section.

\item {\bf Experiments compute resources}
    \item[] Question: For each experiment, does the paper provide sufficient information on the computer resources (type of compute workers, memory, time of execution) needed to reproduce the experiments?
    \item[] Answer: \answerYes{} 
    \item[] Justification: Hardware, per-epoch training time, and the total cost of the reported training run are reported in \cref{subsec:resources} for CT reconstruction and in \cref{sec:matrixresources} for robust matrix completion. The full research project, including preliminary experiments at smaller depths and architecture searches, required several times the compute of the reported runs, as noted in the ``Previous experiments'' paragraph.
    
\item {\bf Code of ethics}
    \item[] Question: Does the research conducted in the paper conform, in every respect, with the NeurIPS Code of Ethics \url{https://neurips.cc/public/EthicsGuidelines}?
    \item[] Answer: \answerYes{} 
    \item[] Justification: We have reviewed the NeurIPS Code of Ethics and confirm that our research conforms with the guidelines in every respect.

\item {\bf Broader impacts}
    \item[] Question: Does the paper discuss both potential positive societal impacts and negative societal impacts of the work performed?
    \item[] Answer: \answerNA{} 
    \item[] Justification: The contribution is theoretical and algorithmic: a worst-case bound for inexact alternating minimization, a greedy training scheme, and a PAC-Bayes generalization argument. The CT reconstruction experiments serve only as a numerical instance of the abstract objective \cref{eq:Objective}. The paper makes no clinical claims, releases no diagnostic model, and shows no patient images.
    
\item {\bf Safeguards}
    \item[] Question: Does the paper describe safeguards that have been put in place for responsible release of data or models that have a high risk for misuse (e.g., pre-trained language models, image generators, or scraped datasets)?
    \item[] Answer: \answerNA{} 
    \item[] Justification: Our models are trained to solve large-scale optimization problems and do not have a high risk of misuse.

\item {\bf Licenses for existing assets}
    \item[] Question: Are the creators or original owners of assets (e.g., code, data, models), used in the paper, properly credited and are the license and terms of use explicitly mentioned and properly respected?
    \item[] Answer: \answerYes{} 
    \item[] Justification: We have summarized the third-party assets used with corresponding license and reference in \cref{tab:licenses}. The data for the second experiment on robust matrix completion is entirely synthetic.

\item {\bf New assets}
    \item[] Question: Are new assets introduced in the paper well documented and is the documentation provided alongside the assets?
    \item[] Answer: \answerYes{} 
    \item[] Justification: The supplement contains the source code, the pre-trained dictionary, and a \texttt{README} that documents the directory layout, dependencies, and the commands needed to reproduce each table and figure. The code is also replicated on  a publicly accessible github repository that is linked in th abstract.

\item {\bf Crowdsourcing and research with human subjects}
    \item[] Question: For crowdsourcing experiments and research with human subjects, does the paper include the full text of instructions given to participants and screenshots, if applicable, as well as details about compensation (if any)? 
    \item[] Answer: \answerNA{} 
    \item[] Justification: Our research does not involve crowdsourcing nor research with human subjects.

\item {\bf Institutional review board (IRB) approvals or equivalent for research with human subjects}
    \item[] Question: Does the paper describe potential risks incurred by study participants, whether such risks were disclosed to the subjects, and whether Institutional Review Board (IRB) approvals (or an equivalent approval/review based on the requirements of your country or institution) were obtained?
    \item[] Answer: \answerNA{} 
    \item[] Justification: No crowdsourcing or research with human subjects was part of this research so there is no requirement for an approval by an ethics board.

\item {\bf Declaration of LLM usage}
    \item[] Question: Does the paper describe the usage of LLMs if it is an important, original, or non-standard component of the core methods in this research? Note that if the LLM is used only for writing, editing, or formatting purposes and does \emph{not} impact the core methodology, scientific rigor, or originality of the research, declaration is not required.
    \item[] Answer: \answerNA{} 
    \item[] Justification: In the paper we used LLMs only for standard components and editing/writing as well as a verifier for the proofs and theorem statements that are included in the paper.

\end{enumerate}